\documentclass[article]{elsarticle}

\usepackage{lineno,hyperref}
\usepackage{import}
\usepackage{placeins}
\usepackage{subcaption}
\modulolinenumbers[5]

\journal{arXiv}

\usepackage{a4}
\usepackage{amssymb}
\usepackage{amsmath}
\usepackage{amsthm}

\usepackage{bm}
\usepackage{color}

\newcommand{\R}{\mathbb R}

\newcommand{\dz}{\mathop{\mathrm{d}z}}
\newcommand{\dx}{\mathop{\mathrm{d}x}}
\newcommand{\dy}{\mathop{\mathrm{d}y}}
\newcommand{\ds}{\mathop{\mathrm{d}s}}
\newcommand{\dOmega}{\mathop{\mathrm{d}\Omega}}

\newcommand{\dr}{\mathop{\mathrm{d}r}}
\newcommand{\dth}{\mathop{\mathrm{d}\theta}}

\newtheorem{remark}{Remark}

\begin{document}

\begin{frontmatter}

\title{A Partitioned Finite Element Method for power-preserving
    discretization of open systems of conservation laws}

\author{Fl\'avio Luiz Cardoso-Ribeiro}
\address{Instituto Tecnol\'ogico de Aeron\'autica, Brazil}

\author{Denis Matignon}
\address{ISAE-SUPAERO, Universit\'e de Toulouse, France}

\author{Laurent Lef\`evre}
\address{{ Univ. Grenoble Alpes, Grenoble INP, LCIS, F-26000, Valence}}




\begin{abstract}
{ This paper presents a structure-preserving spatial discretization method for distributed parameter port-Hamiltonian systems. The class of considered systems are hyperbolic systems of two conservation laws in arbitrary spatial dimension and geometries. For these systems, a partioned finite element method is derived, based on the integration by parts of one of the two conservation laws written in  weak form. The nonlinear 1D Shallow Water Equation (SWE) is first considered as a motivation example. Then the method is investigated on the example of the nonlinear 2D SWE. Complete derivation of the reduced finite-dimensional port-Hamiltonian system is provided and numerical experiments are performed. Extensions to curvilinear (polar) coordinate systems, space-varying coefficients  and higher-order port-Hamiltonian systems (Euler-Bernoulli beam equation) are provided.}
\end{abstract}


\end{frontmatter}


\section{Introduction}

The port-Hamiltonian formalism has been proven to be a powerful tool for the modeling and control of complex multiphysics systems. In many cases, spatio-temporal dynamics must be considered and  infinite-dimensional port-Hamiltonian models are needed. Classical academic examples such as the transmission line, the shallow water or the beam equations have been investigated in the port-Hamiltonian framework \citep{duindam2009modeling}.

Besides 2D and 3D problems have been recently considered \citep{Wu2015b,Vu2016,Trenchant2017}.  In many of these examples, e.g. those arising from mechanics, systems of two  balance equations are considered such as mass and momentum or volume and momentum balance equations.

In order to simulate and design control laws, obtaining a finite-dimensional approximation which preserves the port-Hamiltonian structure of the original system can be advantageous. It may serve as a design guide such as in Control by Interconnection (CbI) or in Interconnection and Damping Assignment Passivity Based Control (IDA-PBC). Besides, preserving the underlying Dirac interconnection structure results in energy conservation properties and associated dynamical properties (e.g. stability, controllability, etc.).

Mixed finite element methods were introduced a long time ago to perform structure-preserving spatial discretization of the Maxwell field equations \cite{bossavit1988whitney,bossavit1998computational}. An extension to open port-Hamiltonian systems was presented in \cite{Golo2004} where the authors proposed a mixed finite element structure-preserving spatial discretization for 1D hyperbolic systems of conservation laws, making use of disctinct  low-order Whitney bases functions to approximate respectively the energy and co-energy variables. These ideas were applied later for the discretization of a parabolic diffusion problem related to pressure-swing-adsorption columns \cite{baaiu2009structure}, piezo-electric beams \cite{CardosoRibeiro2016} and connected to finite volume and staggered grid finite difference methods for 1D problems \cite{Kotyczka2016}. They were also generalized for 2D systems in \cite{Wu2015b} where they are applied to  a vibro-accoustic systems and in \cite{Trenchant2017} where connection was made with 2D finite difference staggered grids schemes. Finally, these structure-preserving mixed finite element methods, applied to the spatial discretization of general port-Hamiltonian systems with boundary energy flows, were stated in a geometry independent form making use of discrete exterior calculus results \cite{seslija2014explicit}. 

In these previous works, the central idea was to define different discretization bases for the energy and co-energy variables such that the strong-form of the equations was exactly satisfied in the corresponding spanned finite-dimensional approximation spaces. This idea was extended to geometric pseudo-spectral methods using conjugated high-order polynomial bases (\cite{Moulla2012}) or Bessel (\cite{Minh2017}) bases functions, globally defined on the whole spatial domain. Still the same idea has been considered in \cite{Farle2014,Farle2013,Farle2014b} for the 1D transmission line equation and for the Maxwell equations. In these latter works, one of the balance law is  kept in strong-form (with exact spatial derivation) while the other one  is considered in the weak sense only. As it was noticed in \cite{Hiemstra2014}, defining these compatible spaces - with power-conjugated approximation bases for the energy and co-energy variables - is (relatively) straightforward for 1D systems, but seems to be cumbersome for higher spatial dimensions or higher order methods. Indeed, the kernel of the exterior derivatives in $N$-D dimensional domains is not anymore trivial and the discretization of the trace operator on boundaries with non trivial (i.e. non rectangular) geometries often leads to dimensionality problems. As suggested in \cite{Kotyczka2017} the discretization of the weak formulation of the considered port-Hamiltonian system may be a practical solution to deal with these higher dimensional problems or more complex geometries. We propose in this paper to follow this approach but to perform integration by parts - which was used in \cite{Kotyczka2017} to get the weak formulation - on one of the two balance equations only, defining in this way a  partitioned mixed finite element method. Doing so, the discretization in the chosen bases for the energy and co-energy variables (and the associated test functions) directly leads to a finite-dimensional Dirac interconnection structure and no further projection is required to get a finite-dimensional port-Hamiltonian system equations with reversible causality. Besides boundary conditions are naturally handled, even in the case of higher-order finite element bases. Finally, the use of this weak-form formulation enables the use of standard finite-element software to perform the proposed structure-preserving discretization and consequently paves the way for further applications to more involved higher-order problems with complex geometries.

 This paper starts with a motivation example detailed in Section \ref{sec:introexample}: a structure-preserving spatial discretization for the one-dimensional (1D) Shallow Water Equations (SWE).  In Section \ref{sec:GeneralSetting}, the approach is generalized in a general setting, where the initial model is stated in a coordinate-free form, independent of the specific geometry. The covariant formulation for systems of two conservation laws with boundary energy flows is presented, as well as its weak form. In Section \ref{sec:2DPFEM} the proposed partitioned finite element method (PFEM) is applied to the 2D SWE example. Numerical experiments are presented in Section \ref{sec:numericalresults}. Finally, extensions of the method to curvilinear (polar) coordinate systems, space-varying coefficients and higher-order port-Hamiltonian systems (Euler-Bernouilli beam equations) are provided in Section \ref{sec:extensions}.   The paper ends with conclusions and open questions which are discussed in Section~\ref{sec:conclusions}.

%
%
%

\section{An introductory example}
\label{sec:introexample}
%
%
%
{ The aim of this section is to present the general idea of this paper - partial integration by parts of the weak form for systems of two conservation laws and structure-preserving projections in finite-dimensional approximation spaces - applied on a simple one-dimensional example, namely the 1D SWE written in the port-Hamiltonian formulation. First, the port-Hamiltonian formulation for these equations is recalled (subsection \ref{subsec:phSWE}). Then the partial integration by parts idea is performed on the weak form for this 1D SWE example and a structure-preserving finite element spatial discretization method is applied to obtain the finite-dimensional port-Hamiltonian system (pHs) (subsection \ref{subsec:wfSWE}).
This general idea differs from previous works (as \cite{Golo2004} and \cite{Moulla2012}) where the central idea was to define different discretization spaces for the energy and co-energy variables such that the strong form of the equations were exactly satisfied in these finite-dimensional spaces. Instead, we use a weak-form representation for the equations, where only one of the conservation laws is integrated by parts. This partitioned approach naturally leads to a skew-symmetric interconnection matrix between the energy and co-energy variables.  Furthermore, the use of weak form enables the use of classical finite-element methods to perform the discretization.
}

{ 
\subsection{Port-Hamiltonian strong formulation for the 1D SWE}\label{subsec:phSWE}}
The Shallow Water Equations are a set of partial differential equations that can be used to represent an incompressible fluid with free-surface motion. These equations are typically used to model fluid motion in water channels \cite{Hamroun2010}, wave propagations in oceans and lakes, and sloshing in fluid tanks \cite{AlemiArdakani2016,Cardoso-Ribeiro2017}. { When one considers the frictionless flow in a horizontal channel with uniform rectangular cross-section, the one-dimensional mass and momentum balance equations may be written as:}
\begin{equation}
	\label{eq:SWE}
	\begin{split}
		\frac{\partial}{\partial t}h = - \frac{\partial}{\partial z}\left(h u\right) \,,\\
		\frac{\partial}{\partial t}u = - \frac{\partial}{\partial z}\left(\frac{u^2}{2} + gh\right) \,,
	\end{split}
\end{equation}
where $h(z,t)$ is the fluid height, $u(z,t)$, { the fluid average velocity in a cross-section}, $z$, the spatial coordinate, $t$, the time and $g$, the gravitational acceleration. 

The total energy of the system inside the domain $[0,L]$ is given by the sum of kinetic and potential (gravitational) energy:
\begin{equation}
  \label{eq:SWE1H}
	H = \frac{1}{2}\int_{[0,L]} {  \left( \rho b \, h u^2 + \rho b g \, h^2  \right)} \dz \,,
\end{equation}
where $b$ is the width of the water channel (or tank) rectangular cross-section and $\rho$, the water density.  Defining the energy-variables $q(z,t) := bh(z,t)$ and $p(z,t) := \rho u (z,t)$, the system Hamiltonian (total energy) is given by:
 \begin{equation}
	\label{eq:waveequation_Ham}
	H\left( q(z,t), p(z,t)\right) = \frac{1}{2}\int_{[0,L]} \left( \frac{q p^2}{\rho} + \, \frac{\rho g}{b} q^2 \right) \dz \,.
 \end{equation}
 Using these newly defined variables, \eqref{eq:SWE} can be rewritten as:
 \begin{equation}
	\label{eq:1Dstrongform}
	\begin{split}
		\dot{q}(z,t) & = - \frac{\partial}{\partial z} e_p (z,t) \,,\\
		\dot{p}(z,t) & = - \frac{\partial}{\partial z} e_q (z,t)\,,
	\end{split}
 \end{equation}
{where $e_q(z,t)$ and $e_p(z,t)$ are the co-energy variables (respectively, the total pressure and the water flow) which are defined as the variational derivatives of the Hamiltonian with respect to $q(z,t)$ and $p(z,t)$:}
 \begin{equation}
	\label{eq:1Dcoenergyvariables}
	\begin{split}
		e_q & = \frac{\delta H}{\delta q} = \frac{p^2}{2\rho} + \frac{\rho g}{b} q = \rho \left( \frac{u^2}{2} + g h \right) \,, \\
		e_p & = \frac{\delta H}{\delta p} = \frac{q p}{\rho}  = b h u \,.
	\end{split}
 \end{equation}
%

%
%
{ From the above definitions of energy and co-energy variables, using the SWE written in the canonical Hamiltonian form (\ref{eq:1Dstrongform}) and Stokes theorem, one obtains for the power balance equation:
\begin{equation}
	\label{eq:Hdot1D}
	\begin{split}
	\dot{H} (t) & = \int_{[0, L]}{\left( e_q(z,t) \dot{q}(z,t) + e_p(z,t) \dot{p}(z,t) \right) \dz} \,, \\	
	& = - \int_{[0,L]}{ \frac{\partial  }{\partial z}  \left( e_q(z,t) e_p(z,t) \right) \dz} \,, \\
	& = - \int_{\partial [0,L]}{e_q(z,t) e_p(z,t)} \,, \\
	& = \pmb{u}_\partial^T \pmb{y}_\partial \,,
	\end{split}
\end{equation}
%
where boundary port input variables, $\pmb{u}_\partial$, are defined as the values of the co-energy variables evaluated at the spatial domain boundary:
\begin{equation}
	\label{eq:input1d}
	\pmb{u}_\partial := \left[\begin{matrix} e_p(0,t) \\ e_p(L,t) \end{matrix} \right] \,,
\end{equation}
while the power-conjugated boundary output variables are defined as:
\begin{equation}
	\label{eq:output1d}
	\pmb{y}_\partial = \left[ \begin{matrix} e_q(0,t) \\ -e_q(L,t) \end{matrix} \right] \,.
\end{equation}
The power balance equation (\ref{eq:Hdot1D}) defines a natural pairing or bilinear form
\begin{equation}
	\label{eq:powerpairing}
 	\begin{split}
	\left< \cdot \left| \cdot \right. \right>: & \; \mathcal{B} \rightarrow \R  \,, \\
	& (e,f) \mapsto \left< e \left| f \right. \right> := \int_{[0, L]}{\left( e_q(z,t) f_q(z,t) + e_p(z,t) f_p(z,t) \right) \dz} + \pmb{u}_\partial^T \pmb{y}_\partial \,,
	\end{split}
\end{equation}
where the Bond space $\mathcal{B}:=\mathcal{E} \times \mathcal{F}$ is defined as the product of the effort real vector space 
\begin{equation}
	\label{eq:effort_space}
	\mathcal{E}:=\left\{ \left. e:=\left[ e_q \; e_p \; \pmb{e}_\partial \right]^T \right| e_q,e_p \in H^1\left( 0,L \right) \; ; \; \pmb{e}_\partial \in \R^2 \right\} \,,
\end{equation}
and its dual flow real vector space
\begin{equation}
	\label{eq:flow_space}
	\mathcal{F}:=\left\{ \left. f:=\left[ f_q \; f_p \; \pmb{f}_\partial \right]^T \right| f_q,f_p \in L^2\left( 0,L \right) \; ; \; \pmb{f}_\partial \in \R^2 \right\} \,,
\end{equation}
with $H^1\left( 0,L \right)$ and   $L^2\left( 0,L \right)$ denoting respectively the Sobolev space of functions with square integrable derivatives on $[0,L]$ and the usual Lebesgue space of square integrable functions on $[0,L]$. Using the bilinear form (\ref{eq:powerpairing}), the power balance equation (\ref{eq:Hdot1D}) simply reads
\begin{equation}
	\label{eq:powerbalance2}
 	\left< \left[\begin{matrix} e_q(t,\cdot) \\ e_p(t,\cdot) \\  \pmb{e}_\partial(t) \end{matrix} \right] \left| \left[\begin{matrix} \dot{q}(t,\cdot) \\ \dot{p}(t,\cdot) \\  \pmb{f}_\partial(t) \end{matrix} \right] \right. \right> = 0 \,,
\end{equation}
with $\pmb{e}_\partial^T(t)=[e_q(0,t) \; e_q(L,T)]$ and $\pmb{f}_\partial^T(t)=[e_p(0,t) \; e_p(L,T)]$ (or the reverse). Besides, the pairing (\ref{eq:powerpairing}) may be symmetrized to obtain the associated inner product:
\begin{equation}
	\label{eq:innerproduct}
 	\begin{split}
	\ll \cdot \left| \cdot \right. \gg : & \; \mathcal{B} \times \mathcal{B} \rightarrow \R  \\
	& \left( (e_1,f_1),(e_2,f_2) \right) \mapsto \ll (e_1,f_1),(e_2,f_2) \gg := \frac{1}{2} \left( \left< \left. e_1  \right| f_2 \right> + \left< \left. e_2  \right| f_1 \right>\right) \,.
	\end{split}
\end{equation}
It may be shown that the Hamiltonian formulation (\ref{eq:1Dstrongform}) for the SWE, together with the boundary conditions (\ref{eq:input1d}) and output (\ref{eq:output1d}) may be equivalently implicitely defined as 
\begin{equation}
	\label{eq:implSDS1D}
\left( \left(\frac{\delta H}{\delta q} ,\frac{\delta H}{\delta p},\pmb{u}_\partial \right),  \left( -\frac{dq}{dt} ,- \frac{dp}{dt},\pmb{y}_\partial \right) \right) \in \mathcal{D} \,,
\end{equation}
where $\mathcal{D}\subset \mathcal{B}$ is the linear subspace which is maximally isotropic (i.e. $\mathcal{D}= \mathcal{D}^\perp$) with respect to the inner product (\ref{eq:innerproduct}) \cite{Moulla2012}. {In that sense, the natural pairing (\ref{eq:powerpairing}) fully describes the geometric structure of the port-Hamiltonian system~(\ref{eq:1Dstrongform}) with boundary values~(\ref{eq:input1d}).} Therefore, in this paper, structure-preserving (or symplectic) spatial discretization will be understood as approximations (projections) which preserve this power form (\ref{eq:powerpairing}). Symplecticity in that sense not only implies preservation of the power balance (\ref{eq:Hdot1D}) or (\ref{eq:powerbalance2}) (i.e. isotropy) but preservation of the whole geometric structure of the system (e.g. the Poisson structure in the example of closed systems or the Dirac structure in the case of open systems with time varying boundary conditions) \cite{Moulla2012,Kotyczka2017}.
}
%

Note that the particular input and output port variables chosen here above in (equations \eqref{eq:input1d} and \eqref{eq:output1d}) is only one among other possible choices. A description of all the possible choices of input/output variables which lead to well-posed problems (in the {\em linear} case) is described in \cite{Gorrec2005}. 

%
\subsection{Partitioned weak form and structure-preserving discretization for the 1D SWE}\label{subsec:wfSWE}
{ We will now introduce a weak formulation for the 1D SWE and then perform integration by parts on the mass balance equation. Let $v_q(z) \in H^1\left( 0,L \right)$  and $v_p(z) \in L^2\left( 0,L \right)$ denote any arbitrary test functions, we may obtain from the strong formulation \eqref{eq:1Dstrongform} the following weak form:
\begin{equation}
	\begin{split}
		\int_{[0,L]}{v_q(z) \dot{q}(z,t) \dz\,} & = -\int_{[0,L]} v_q(z) \frac{\partial}{\partial z} e_p (z,t) \dz \,,\\ 		\int_{[0,L]}{v_p(z) \dot{p}(z,t) \dz\,} & = - \int_{[0,L]} v_p(z) \frac{\partial}{\partial z} e_q (z,t) \dz \,.
	\end{split}
 \end{equation}
Integrating by part the {\em mass balance equation} only, we get the following partitioned weak form:
\begin{equation}
\label{eq:weakform1D}
	\begin{split}
		\int_{[0,L]}{v_q(z) \dot{q}(z) \dz\,} & = \int_{[0,L]} e_p (z,t) \frac{\partial}{\partial z} v_q(z)   \dz\, - v_q(L) e_p(L,t) + v_q(0) e_p(0,t)  \,, \\
		\int_{[0,L]}{v_p(z) \dot{p}(z) \dz\,} & = - \int_{[0,L]} v_p(z) \frac{\partial}{\partial z} e_q (z,t) \dz \,.
	\end{split}
\end{equation}
\begin{remark}
In the specific case where $v_q(z) = 1$ and $v_p(z) = 1$, we get:
	\begin{equation}
	\begin{split}
		\int_{[0,L]}{ \dot{q}(t,z) \dz\,} & =  e_p(0,t) -  e_p(L,t) \,, \\
		\int_{[0,L]}{ \dot{p}(t,z) \dz\,} & = e_q(0,t) - e_q(L,t) \,,
	\end{split}
	\end{equation}
	which shows that the two conservation laws for the total mass and the total momentum in the spatial domain $[0,L]$ are preserved in the weak formulation. When $v_q = e_q(z,t)$ and $v_p = e_p(z,t)$ are chosen, one gets:
\begin{equation}
\label{eq:weakform1DHconservation}
	\begin{split}
		\int_{[0,L]}{e_q(z,t) \dot{q}(z) \dz\,} & = \int_{[0,L]} e_p (z,t) \frac{\partial}{\partial z} e_q(z)   \dz\, - e_q(L,t) e_p(L,t) + e_q(0,t) e_p(0,t) \,, \\
		\int_{[0,L]}{e_p(z,t) \dot{p}(z) \dz\,} & = - \int_{[0,L]} e_p(z) \frac{\partial}{\partial z} e_q (z,t) \dz \,.
	\end{split}
\end{equation}
Therefore the power-balance equation \eqref{eq:Hdot1D} reads:
\begin{equation}
	\label{eq:powerbalance2d}
		\dot{H}   = \int_{[0,L]}{(e_q(z,t) \dot{q}(z) + e_p(z,t) \dot{p}(z)) \dz} =  - e_q(L,t) e_p(L,t)  + e_q(0,t) e_p(0,t) \,,
\end{equation}
which shows that the power balance is also preserved in the weak formulation.
\end{remark}
}

{ We will now project the partitioned weak formulation (\ref{eq:weakform1D}) into finite-dimensional approximation spaces chosen in such a way as to preserve the total mass and momentum conservation laws, the power balance equation and the underlying Dirac structure of the original port-Hamiltonian model \eqref{eq:1Dstrongform}. Unlike in \cite{Moulla2012,Kotyczka2017} where different approximation bases are chosen for the energy and co-energy variables, we obtain the mass, momentum, power and structure-preservation by the selection of different approximation bases for the mass and momentum densities.  This ``partitioned'' choice for the approximation bases lead us to square skew-symmetric interconnection matrices. 

Let us approximate the energy variables $q(z,t)$ and $p(z,t)$ as 
\begin{equation}
	\label{eq:approx}
	\begin{split}
	q(z,t) \approx q^{ap}(z,t) & := \sum_{i=1}^{N_q} \phi^i_q(z) q^i(t) = \pmb{\phi}^T_q(z) \pmb{q}(t) \,, \\
	p(z,t) \approx p^{ap}(z,t) & := \sum_{j=1}^{N_p} \phi^j_p(z) p^j(t) = \pmb{\phi}^T_p(z) \pmb{p}(t) \,,
\end{split}
\end{equation}
where  $\phi^i_q(z)\, , \; i\in\{ 1,\ldots, N_q\}$ are the chosen approximation basis functions in $H^1(0,L)$, $\phi^i_p(z)\, , \; i\in\{ 1,\ldots, N_p\}$ the chosen approximation basis functions in $L^2(0,L)$, while $\pmb{q}(t)$ and $\pmb{p}(t)$ are the approximation coordinates for $q^{ap}(z,t)$ and $p^{ap}(z,t)$ in the approximation bases $\pmb{\phi}_q(z)$ and $\pmb{\phi}_p(z)$. The test functions $v_q(z)$ and $v_p(z)$ are approximated in the same bases as $q(z,t)$ and $p(z,t)$, respectively.}
%
%
From the substitution of the approximated variables \eqref{eq:approx}  in the weak form \eqref{eq:weakform1D}, the following finite-dimensional equations are obtained:
\begin{equation}
\label{eq:weakform1Dwithv}
	\begin{split}
		\pmb{v_q}^T\left[ \int_{[0,L]}{\pmb{\phi}_q(z) \pmb{\phi}^T_q(z) \dz} \right] \dot{\pmb{q}}(t)  = & \; \pmb{v_q}^T \left[  \int_{[0,L]}{\frac{d \pmb{\phi}_{q}}{dz}(z) \pmb{\phi}_{p}^T(z) \dz} \right] \pmb{e}_p (t) \\
		                        & + \pmb{v}_q^T \pmb{\phi}_q(0) e_p^{ap}(0,t) - \pmb{v}_q^T \pmb{\phi}_q(L) e_p^{ap}(L,t) \,, \\
		\pmb{v_p}^T \left[ \int_{[0,L]}{\pmb{\phi}_p(z) \pmb{\phi}^T_p(z) \dz} \right] \dot{\pmb{p}}(t)   =& - \pmb{v_p}^T \left[ \int_{[0,L]}{\pmb{\phi}_p(z) \frac{d \pmb{\phi}^T_{q}}{dz}(z) \dz} \right] \pmb{e}_q (t) \,,
	\end{split}
\end{equation}
{where the $H^1(0,L)$ effort functions $e_q(\cdot,t)$ have been approximated in the  $\{ \phi^i_q(z) \}$ basis, while  the $L^2(0,L)$ effort functions $e_p(\cdot,t)$ have been approximated in the  $\{ \phi^i_p(z) \}$ basis. Similarly, the flow functions $\dot{q}(\cdot,t)$ and $\dot{p}(\cdot,t)$ have been approximated respectively in the  $\{ \phi^i_q(z) \}$ and $\{ \phi^i_p(z) \}$ bases. Since these equations should remain valid for any choices of test functions coordinates $\pmb{v}_1$ and $\pmb{v}_2$, one gets:
\begin{equation}
	\label{eq:finitedim_waveeq1d}
	\begin{split}
		{M_q} \dot{\pmb{q}}(t) & = {D} \pmb{e}_p (t) + B \left[\begin{matrix}e_p(0,t) \\e_p(L,t) \end{matrix} \right] \,, \\
		{M_p} \dot{\pmb{p}}(t)  & = {-D^T} \pmb{e}_q (t)\,,
	\end{split}
\end{equation}
where $M_q$ and $M_p$ are square mass matrices (of size $N_q \times N_q$ and $N_p \times N_p$, respectively) defined as 
 \begin{equation}
	\label{eq:MqMpDef}
		{M_q} :=  \int_{[0,L]}{\pmb{\phi}_q(z) \pmb{\phi}^T_q(z) \dz} \,, \hspace{2cm} 
		{M_p} :=  \int_{[0,L]}{\pmb{\phi}_p(z) \pmb{\phi}^T_p(z) \dz} \,.
\end{equation}
Matrix $D$ is of size $N_q \times N_p$, and defined as 
 \begin{equation}
	\label{eq:DDef}
		{D}  := \int_{[0,L]}{\frac{d \pmb{\phi}_{q}}{dz}(z) \pmb{\phi}^T_{p}(z) \dz} \,,
\end{equation}
and $B:=\left[\begin{matrix}\pmb{\phi}_q(0) & - \pmb{\phi}_q(L) \end{matrix} \right]$ is an $N_q \times 2$ matrix.  Using the input-output conjugated boundary port variables as defined in \eqref{eq:input1d} and \eqref{eq:output1d}, since the boundary values of $e_q^{ap}(z,t)$ may be written
\begin{equation}
	\left[
		\begin{matrix}
			e_q^{ap}(0,t) \\ -e_q^{ap}(L,t)
		\end{matrix}
	\right]
	=
	\left[\begin{matrix}\pmb{\phi}^T_q(0) \\ - \pmb{\phi}^T_q(L) \end{matrix} \right] \pmb{e}_q = B^T \pmb{e}_q
\end{equation}
the approximation \eqref{eq:finitedim_waveeq1d} may be written using the following finite-dimensional Dirac structure representation:
\begin{equation}
	\label{eq:finitedim_waveeq}
	\begin{split}
		\left[ \begin{matrix} M_q & 0 \\ 0 & M_p \end{matrix}\right] \left[ \begin{matrix} {\pmb{f}}_q(t) \\ {\pmb{f}}_p(t) \end{matrix}\right] & = \left[ \begin{matrix} 0 & D \\ -D^T & 0 \end{matrix}\right] \left[ \begin{matrix} {\pmb{e}}_q(t) \\ {\pmb{e}}_p(t) \end{matrix}\right] + \left[ \begin{matrix} B \\ 0 \end{matrix}\right] \left[\begin{matrix}e_p(0,t) \\e_p(L,t) \end{matrix} \right] \,,\\
		\left[
		\begin{matrix}
			e_q^{ap}(0,t) \\-e_q^{ap}(L,t)
		\end{matrix}
	\right] & = \left[ \begin{matrix} B^T & 0 \end{matrix}\right] \left[ \begin{matrix} {\pmb{e}}_q(t) \\ {\pmb{e}}_p(t) \end{matrix}\right] \,,
	\end{split}
\end{equation}
where ${\pmb{f}}_q(t)$ and ${\pmb{f}}_p(t)$ denote the vector coordinates for the flow approximations respectively in the  in the $\phi^i_q(z)$ and $\phi^j_p(z)$ approximation bases, that is $f_q^{ap}(z,t) = \pmb{\phi}_q^T(z) \pmb{f}_q(t)$ and $f_p^{ap}(z,t) = \pmb{\phi}_p^T(z) \pmb{f}_p(t)$.

We obtained the finite-dimensional Dirac structure representation \eqref{eq:finitedim_waveeq} from the projection of the Stokes-Dirac structure on a the chosen approximation spaces. We will now derive the corresponding approximation of the port-Hamiltonian system dynamics \eqref{eq:implSDS1D} by restricting the Hamiltonian functional to the same approximation spaces. From the definition of co-energy variables as variational derivatives of the Hamiltonian with respect to $q$ and $p$, we get:
\begin{equation}
	\dot{H}(t) = \int_{[0,L]} \left( e_q(z,t) \dot{q}(z,t)  + e_p(z,t) \dot{p}(z,t)  \right)) \dz \,.
\end{equation}
Using the approximations \eqref{eq:approx} for the energy and co-energy variables, this power balance may be approximated as:
\begin{equation}
\label{eq:continousHd}
\begin{split}
	\dot{H_{d}}(t) & :=\int_{[0,L]} \left( 
e_q^{ap}(z,t) \dot{q}^{ap}(z,t)  + e_p^{ap}(z,t) \dot{p}^{ap}(z,t)  
\right) \dz \,,
\\
 & = \pmb{e}_q^T(t) M_q \dot{\pmb{q}}(t) + \pmb{e}_p^T(t) M_p \dot{\pmb{p}}(t) \,, \\ &=
\pmb{e}_q^T(t) \dot{\tilde{\pmb{q}}}(t) + \pmb{e}_p^T(t) \dot{\tilde{\pmb{p}}}(t) \,,
\end{split}
\end{equation}
where new energy variables
\begin{equation}
	\begin{split}
	\tilde{\pmb{q}}(t) & := M_q \pmb{q}(t) \,, \\
	\tilde{\pmb{p}}(t) & := M_p \pmb{p}(t) \,,
	\end{split}
\end{equation}
have been defined. Therefore, in order to write the power balance \eqref{eq:continousHd} as the total time derivative of the discrete Hamiltonian written as a function of the finite-dimensional vector coordinates for the energy variables, the following reduced co-energy variables must be defined: 
\begin{equation}
	\label{eq:0DSWEconstrel1}
	\begin{split}
		{\pmb{e}}_q(t) & = \frac{\partial \tilde{H}}{\partial \tilde{\pmb{q}}} \,,\\
		{\pmb{e}}_p(t) & = \frac{\partial \tilde{H}}{\partial \tilde{\pmb{p}}}\,,
	\end{split}
\end{equation}
where
\begin{equation}
\label{eq:0DSWEconstrel2}
	\tilde{H}(\tilde{\pmb{q}},\tilde{\pmb{p}}):= H\left(  \pmb{\phi}^T_q(z) M_q^{-1} \tilde{\pmb{q}}(t), \pmb{\phi}^T_p(z) M_p^{-1} \tilde{\pmb{p}}(t) \right) \,.
\end{equation}

We obtain the finite-dimensional pHs formulation for the proposed structure-preserving reduction scheme by combining equations \eqref{eq:finitedim_waveeq} (for the linear finite-dimensional Dirac interconnection structure) and the nonlinear constitutive equations \eqref{eq:0DSWEconstrel1} and \eqref{eq:0DSWEconstrel2}:
\begin{equation}
	\label{eq:0DSWEPCH}
	\begin{split}
		\left[ \begin{matrix} {\dot{\tilde{\pmb{q}}}}(t) \\ {\dot{\tilde{\pmb{p}}}}(t) \end{matrix}\right] & = \left[ \begin{matrix} 0 & D \\ -D^T & 0 \end{matrix}\right] \left[ \begin{matrix} \frac{\partial \tilde{H}}{\partial \tilde{\pmb{q}}} \\ \frac{\partial \tilde{H}}{\partial \tilde{\pmb{p}}} \end{matrix}\right] + \left[ \begin{matrix} B \\ 0 \end{matrix}\right] \left[\begin{matrix}e_p(0,t) \\e_p(L,t) \end{matrix} \right] \,,\\
		\left[
		\begin{matrix}
			e_q^{ap}(0,t) \\- e_q^{ap}(L,t)
		\end{matrix}
	\right] & = \left[ \begin{matrix} B^T & 0 \end{matrix}\right] \left[ \begin{matrix} \left[ \begin{matrix} \frac{\partial \tilde{H}}{\partial \tilde{\pmb{q}}} \\ \frac{\partial \tilde{H}}{\partial \tilde{\pmb{p}}} \end{matrix}\right] \end{matrix}\right] \,.
	\end{split}
\end{equation}

The state space model \eqref{eq:0DSWEPCH} for the reduced dynamics exhibits the usual general pHs form for finite-dimensional systems \cite{VanderSchaft2014}. It is a generic formulation which may be obtained for any 1D systems of conservation laws written in the port-Hamiltonian formulation. In the 1D SWE example, the Hamiltonian function is neither quadratic, nor separable. Nevertheless, an explicit form may be obtained for the constitutive equations \eqref{eq:0DSWEconstrel1}. Since the Hamiltonian function restricted to the approximation spaces for $q$ and $p$ reads:
\begin{equation}
\label{eq:0DSWEconstrel3}
	\tilde{H}(\tilde{\pmb{q}},\tilde{\pmb{p}}):=\frac{1}{2}\int_{[0,L]} \left( \frac{\pmb{\phi}^T_q(z) M_q^{-1} \tilde{\pmb{q}} \left( \pmb{\phi}^T_p(z) M_p^{-1} \tilde{\pmb{p}} \right)^2}{\rho} + \, \frac{\rho g}{b} \left( \pmb{\phi}^T_q(z) M_q^{-1} \tilde{\pmb{q}} \right)^2 \right) \dz \,.
\end{equation}
One obtains for the reduced effort variables the expressions:
\begin{equation}
	\label{eq:0DSWEconstrel4}
	\begin{split}
		\frac{\partial \tilde{H}}{\partial \tilde{\pmb{q}}} & = \frac{\rho g}{b}M_q^{-T}\tilde{\pmb{q}}(t)+\left[  \frac{1}{2\rho} \int_{[0,L]}  M_q^{-T}  \pmb{\phi}_q(z) \tilde{\pmb{p}}^T(t)  \pmb{\phi}_p(z)\pmb{\phi}^T_p(z) M_p^{-1}  \dz \right] \tilde{\pmb{p}}(t) \,,\\
		\frac{\partial \tilde{H}}{\partial \tilde{\pmb{p}}} & = \tilde{\pmb{q}}^T(t) \left[  \frac{1}{\rho}\int_{[0,L]}  M_q^{-T}  \pmb{\phi}_q(z)  M_p^{-T} \pmb{\phi}_p(z) \pmb{\phi}_p^T(z)  M_p^{-1}   \dz \right] \tilde{\pmb{p}}(t) \,.
	\end{split}
\end{equation}
Note that both constitutive equations exhibits nonlinear terms. In order to compute them, the following procedured was used. The first equation can be written as:
\begin{equation}
	\frac{\partial \tilde{H}}{\partial \tilde{\pmb{q}}} = M_q^{-T} \underbrace{\left(\frac{\rho g}{b}\tilde{\pmb{q}}(t)+\left[  \frac{1}{2\rho} \int_{[0,L]}   \pmb{\phi}_q(z) \tilde{\pmb{p}}^T(t)  \pmb{\phi}_p(z)\pmb{\phi}^T_p(z)  \dz \right] \tilde{\pmb{p}}(t) \right)}_{\pmb{w}(t)} \,,
\end{equation}
where the components of $\pmb{w}(t)$ can be computed as:
\begin{equation}
	w_i(t) = \frac{\rho g}{b} \tilde{q}_i(t)+ \tilde{\pmb{p}}^T(t) \left(\int_{[0,L]}   {\phi}_{q,i}(z)   \pmb{\phi}_p(z)\pmb{\phi}^T_p(z)  \dz \right) \tilde{\pmb{p}}(t)  \,.
\end{equation}
Note that $\int_{[0,L]}   {\phi}_{q,i}(z)   \pmb{\phi}_p(z)\pmb{\phi}^T_p(z)  \dz$, for $\{i = 1, \dots, N_q\}$  are $N_q$ matrices of dimension $N_p \times N_p$, that can be computed once and remains constant. Similarly, the second constitutive relationship can be written as a function of constant matrices.

}
\begin{remark}
We may deduce from the pHs representation \eqref{eq:0DSWEPCH} that the power balance equation reads:
\begin{equation}
	\dot{H}_d = -e_1(0,t) e_2(0,t) + e_1(L,t) e_2(L,t)  \,.
\end{equation}
Hence, the power balance (and the corresponding power product value) is preserved by the proposed partitioned spatial discretization scheme. In that sense, we call it a structure-preserving or {\em symplectic} scheme. When 
the bases functions satisfy
		\begin{equation}
		   \sum_{i=1}^{N_q}\phi_q^i(z) = \sum_{i=1}^{N_p} \phi_p^i(z) = 1 \, , \; \forall z\in [0,L] \,,
		\end{equation}
the mass and momentum conservation laws are also satisfied in the finite-dimensional approximation spaces: 
	\begin{equation}
		\begin{split}
			{\int_{[0,L]}{ \pmb{\phi}_q(z)^T \dz\,}} \dot{\pmb{q}}(t) = {\int_{[0,L]}{ \dot{\pmb{q}}^{ap}(z,t) \dz}} & = e_p(0,t) -e_p(L,t)  \,, \\
{\int_{[0,L]}{ \pmb{\phi}_p(z)^T \dz\,}} \dot{\pmb{p}}(t) = {\int_{[0,L]}{ \dot{\pmb{p}}^{ap}(z,t) \dz}} & = e_q(0,t) -e_q(L,t)  \,.
		\end{split}
	\end{equation}	
\end{remark}

In the beginning of this section, we motivated this work by the fact that previous work on structure-preserving spatial discretization that relies on exact satisfaction of the strong form of the equations usually lead to difficulties when generalizing to 2D or 3D systems. The following questions arise: can the proposed PFEM method be easily generalized to higher-dimensional problems (2D and 3D)? Does it work with different coordinate systems? What about convergence? We will answer these questions in the following sections.

{

\section{A general setting}
\label{sec:GeneralSetting}
In this section we will define a port-Hamiltonian system of two conservation laws with time-varying port boundary variables in the so-called co-variant form, that is defining energy and co-energy variables as differential forms \cite{VanderSchaft2002}. This will allow us to define, independently from the particular spatial dimension, geometry or coordinate system, the class of problems which can be solved by using the structure-preserving spatial discretization scheme proposed in this paper. In the next subsection (\S~\ref{EDF}), to make the paper self-contained, we give a short introduction to exterior calculus with differential forms and their functional spaces. Then, in \S~\ref{EDFSDS}, we generalize the definition of port-Hamiltonian systems, given in  the previous section for the 1D SWE example, for general systems of two conservation laws, using these differential forms. In \S~\ref{subsec:SDS},  we define de Stokes-Dirac interconnection structures associated to these port-Hamiltonian systems and the corresponding power pairings which will be preserved in the discretization. We also detail the 2D SWE which will be used later in the numerical experiments. Finally, in \S~\ref{EDFWF}, we will give the general weak form which will be used for the structure-preserving discretization scheme presented in section \ref{sec:2DPFEM}. Readers could refer to \cite[chapter 5]{meyer2008}, \cite{flanders1963} or \cite{frankel2011} for a general intuitive introduction to differential forms and their use in physical systems modelling.  Numerical approximations using finite element spaces of differential forms are described in \cite{arnold2010,arnold2013}. 

\subsection{Differential forms and exterior calculus}
\label{EDF}
Let $\Omega$ be an open, bounded and connected $n$-dimensional spatial domain (manifold) with a $(n-1)$-dimensional Lipschitz boundary $\partial \Omega$. 

We denote $\Lambda^{k}\left(\Omega\right)$ the space of smooth differential $k$-forms (i.e. field of alternated $k$-forms defined on the tangent subspace $T_{\xi}\Omega , \, \xi\in\Omega$)  with smooth coefficients in $\Omega$.

Differential $k$-forms are endowed with a product, called {\em exterior
product}, denoted as $\wedge$ , used for instance to define the power as the product of $k$-forms flow variables (time derivatives of energy variables) and $(n-k)$-forms co-energy variables (variational derivatives of the Hamiltonian with respect to the energy variables). The exterior product (also called \emph{wedge product} or
\emph{Grassman product}) is a skew symmetric exterior product such that 
\begin{equation}
\omega^{k}\wedge \omega^{l} = (-1)^{kl} \omega^{l}\wedge \omega^{k} \in\Lambda^{k+l}\left(\Omega\right) \,,
\end{equation}
for all $k$\emph{-form} $\omega^{k}\in\Lambda^{k}\left(\Omega\right)$, $l$\emph{-form} $\omega^{l}\in\Lambda^{l}\left(\Omega\right)$ and $k,l\in \{ 0,\ldots , n\}$ with $k+l\leq n$.

Differential forms are also endowed with a derivation which is called
exterior derivation, denoted by $\textrm{d}$, which generalize the gradient, divergence and curl operators respectively for $0$-forms (i.e. functions), $1$-forms and $2$-forms in $\Omega=\R^3$. The exterior derivative (or \emph{co-boundary map}) is a derivation of degree one
\begin{eqnarray*}
\mbox{d}:\mbox{ }\Lambda^{k}\left(\Omega\right) & \rightarrow & \Lambda^{k+1}\left(\Omega\right) \,,
\end{eqnarray*}
such that
\begin{eqnarray}
\mbox{d}(\omega^k+\eta^k)&=& \mbox{d}\omega^k+\mbox{d}\eta^k \,, \\
\textrm{d}\left(\omega^{k}\wedge\omega^{l}\right)&=&  \textrm{d}\omega^{k}\wedge\omega^{l}+\left(-1\right)^{k} \omega^{k}\wedge\textrm{d}\omega^{l} \label{eq:productrule} \,,\\ 
\textrm{d}\wedge\textrm{d}&=&0 \,,
\end{eqnarray}
for all $\omega^{k},\eta^{k}\in\Lambda^{k}\left(\Omega\right)$ and $\omega^{l}\in\Lambda^{l}\left(\Omega\right)$. The exterior derivative is defined such that the Stokes theorem applies, that is
\begin{equation}
\label{eq:StokesTheorem-1-1} 
\int_{\Omega}\textrm{d}\omega^{n-1}=\int_{\partial\Omega}\textrm{tr} \left( \omega^{n-1}\right) \,,
\end{equation}
where the {\em trace} $\textrm{tr}\left( \omega^{n-1}\right)\in\Lambda^{n-1}\left(\partial \Omega\right)$ denotes the continuous extension of $\omega^{n-1} \in\Lambda^{n-1}\left(\Omega\right)$ to the boundary $\partial \Omega$. The Stokes theorem is the central result which allows us to define, from the natural power pairing, the  associated Stokes-Dirac interconnection structure (see \S~\ref{subsec:SDS}). Symplectic discretization schemes will be defined in order to preserve either the power product or equivalently the Stokes theorem in finite-dimensional approximation spaces. 
In the sequel, we will make use of the duality product (natural power pairing) in $\Lambda^{n}\left(\Omega\right)$ between $k$-forms $\omega^k\in\Lambda^{k}\left(\Omega\right)$ and $(n-k)$-forms $\omega^{n-k}\in \Lambda^{n-k}\left(\Omega\right)$:
\begin{equation}
	\label{eq:nppcv}
\left< \omega^k \right. \left| \omega^{n-k} \right>_\Omega :=\int_{\Omega} \omega^{k}\wedge\omega^{n-k}  \,.
\end{equation}

In this duality product (power pairing), one of the argument (differential form) will denote a flow variable while the other will be an effort variable. These two variables will be related together, on one side by balance equations (conservation laws) which are structural (metric independent) and on the other side by constitutive equations. The constitutive equations are related to the physical properties of the considered material domain (magnetic, visco-elastic, etc.). Therefore, they are metric dependent and will be formulated using the so-called Hodge star operator. For instance, in Section \ref{subsec:wave2d}, the constitutive equations \eqref{eq:HSW2Ddef} or \eqref{eq:eqepSW2Ddef} are derived from the Hamiltonian density which describes, for the 2D SWE, the surfacic potential and kinetic energy densities. These constituve equations relates pressure tensors and water flows with mass and momentum densities. Therefore, metric dependent relations are needed and the Hodge star operator is introduced. The Hodge star operator induces an inner product on the space of differential $k$-forms $L^2\Lambda^k(\Omega)$ on the $n$-dimensional manifold $\Omega$ from the duality product: 
\begin{equation}
\label{eq:Hodgestardef}
\left( \alpha , \beta \right) := \left< \alpha \right. \left| \star{\beta} \right>_\Omega = \left<  {\beta} \right. \left|   \star{\alpha} \right>_\Omega = \left( \beta , \alpha \right) \, , \; \forall \alpha,\beta \in L^2\Lambda^k(\Omega) \,.
\end{equation}
This chosen product is not necessarily the standard inner product for the  $L^2$ norm. For instance, non uniform metric are required when the material domain properties are not uniform in space or isotropic. However, unless otherwise stated, we will consider in this paper the usual standard (metric). For instance, using Cartesian coordinates in some 2-dimensional domain $\Omega$, we will get $\star{1}=dx\wedge dy$, $\star{dx}=dy$ and $\star{dy}=-dx$.

The generalized Stokes' theorem \eqref{eq:StokesTheorem-1-1}, together with the product rule \eqref{eq:productrule}, gives the following integration by parts formula for the power pairing:
\begin{equation}
	\label{eq:nppipp}
\left< \textrm{d}\omega^k \right. \left| \omega^{n-k-1} \right>_\Omega = \left< \textrm{tr } \omega^k \right. \left| \textrm{tr } \omega^{n-k-1} \right>_{\partial \Omega} - (-1)^k \left< \omega^k \right. \left| \textrm{d} \omega^{n-k-1} \right>_\Omega  \,.
\end{equation}
This formula will be helpful, in Section \ref{EDFWF}, to obtain the generalized weak port-Hamiltonian formulation for systems of two conservation laws with boundary energy flows.

When dealing with finite elements approximation spaces, one has to pay some attention to  the regularity of differential forms if we want to apply the previous formula. Readers can refer to \cite[section 4]{arnold2013} for an introduction to differential forms whose coefficients are in Lebesgue spaces $L^p(\Omega)$ and Sobolev spaces $H^m(\Omega)$. We will here point out only two comments on the weak exterior derivative and the trace theorem for differential forms.

Let us consider a $k$-form $\omega^k \in L^2\Lambda^k(\Omega)$, that is a $k$-form $\omega^k \in \Lambda^k(\Omega)$ with coefficients in $L^2(\Omega)$. Then, its {\em weak exterior derivative} $\mathrm{d}\omega^k$ may be defined, using the integration by parts formula \eqref{eq:nppipp}, as the unique $(k+1)$-form (when it exists) such that:
\begin{equation}
	\label{eq:wed}
\left< \textrm{d}\omega^k \right. \left| \mu \right>_\Omega = - (-1)^k \left< \omega^k \right. \left| \textrm{d} \mu \right>_\Omega 
\end{equation}
for all $(n-k-1)$-form $\mu \in C^\infty_c \Lambda^{n-k-1}(\Omega)$, that is a $(n-k-1)$-form $\mu \in \Lambda^{n-k-1}(\Omega)$ with coefficients in the space of $C^\infty$ continuous functions with compact support in $\Omega$. In the sequel of the paper, when dealing with the weak formulation of port-Hamiltonian system, exterior derivatives will be understood in this weak sense.

As far as we are dealing with systems of conservation laws with boundary energy flows, boundary values of the co-energy variables will be necessary to define input/output pairs of power conjugated boundary port variables (see section \ref{EDFSDS} hereafter). Hopefully, we can make use of the so-called {\em trace theorem} which has been extended to Sobolev spaces of differential forms \cite[section 4]{arnold2013}:
\begin{equation}
	\label{eq:tracethm} 
	\omega^k \in H^1\Lambda^k(\Omega) \Rightarrow \textrm{tr }\omega^k \in H^{1/2}\Lambda^k(\partial \Omega)\subset  L^2\Lambda^k(\partial \Omega)\,.
\end{equation}
The injection in equation \eqref{eq:tracethm} is continuous when Lebesgue space $L^2\Lambda^k(\Omega)$ and the Sobolev space $H^1\Lambda^k(\Omega)$ are equipped with the usual inner products. Throughout the paper, we will make use of the following compact notation:
\begin{equation}
\int_{\partial \Omega} \omega^{n-1}:=\int_{\partial\Omega}\textrm{tr} \left( \omega^{n-1}\right)
\end{equation}
for all $\omega^{n-1}$ in $H^m\Lambda^{n-1}(\Omega),\, m\geq 1$.

\subsection{The covariant port-Hamiltonian formulation for systems of two conservation laws with boundary energy flows}
\label{EDFSDS}
We will now extend the port-Hamiltonian formulation which has been presented in Section \ref{subsec:phSWE} only for the 1D SWE example. Hyperbolic systems of two conservation laws will be stated using differential forms, as introduced in the previous section, both for the energy and co-energy variables. Let us consider the two conserved quantities  $\alpha^q \in L^2\Lambda^q(\Omega)$ and ${\alpha^p}\in L^2\Lambda^p(\Omega)$ with $q+p=n+1$ and $n$ denotes the dimension of the open, connected domain $\Omega$ with Lipschitz boundary $\partial \Omega$. These variables $\alpha^q=\alpha^q(\pmb{z},t)$ and $\alpha^p=\alpha^p(\pmb{z},t)$ are vector-valued distributed energy state variables defined for any $\pmb{z}\in\Omega$ ($\pmb{z}$ is the position vector) and time $t\geq 0$. 

Let the Hamiltonian functional $H$ be defined as:
\begin{equation}
	\label{eq:Hdef}
	H(\alpha^q,\alpha^p):=\int_\Omega{\mathcal{H}\left( \alpha^q(\pmb{z},t), \alpha^p(\pmb{z},t),\pmb{z} \right)} \,,
\end{equation}
where $\mathcal{H}$ denotes the Hamiltonian density $n$-form which is assumed to be a smooth function. The variational derivatives of $H$ with respect to $\alpha^q$ and $\alpha^p$ are the unique differential $n-q$ and $n-p$ forms, denoted respectively 
$\delta_{{q}}H$ and $\delta_{{p}}H$, such that:
\begin{equation}
	\label{eq:VDdef}
	H(\alpha^q+\delta \alpha^q, \alpha^p+\delta \alpha^p)=H(\alpha^q,\alpha^p)+\int_\Omega{\left( \delta_{{q}}H \wedge \delta{\alpha^q} + \delta_{{p}}H \wedge \delta{\alpha^p} \right) }+o\left( \delta \alpha^q,\delta \alpha^p \right) \,.
\end{equation}
Therefore, from the Hamiltonian defined in \eqref{eq:Hdef}, we may define the co-energy variables (efforts):
\begin{equation}
	\label{eq:eqpdef}
	\begin{array}{c} {e_q}:=\delta_{{q}}H \,, \\   {e_p}:=\delta_{{p}}H \,, \end{array} 
\end{equation}
with  ${e_q} \in H^1\Lambda^{n-q}(\Omega)$ and ${e_p} \in H^1\Lambda^{n-p}(\Omega)$. The Hamiltonian system of two canonically interacting conservation laws for $\alpha^q$ and $\alpha^p$ may be defined as:
\begin{equation}
	\label{eq:S2CL}
	\left[ \begin{array}{c} \dot{\alpha}^q(\pmb{z},t) \\  \dot{\alpha}^p(\pmb{z},t) \end{array} \right]  = \underbrace{ \left[ \begin{array}{cc} 0 & \mbox{d} \\ (-1)^r  \mbox{d} & 0 \end{array} \right]}_{\mathcal{J}} \left[ \begin{array}{c}  {e_q} (\pmb{z},t) \\ {e_p} (\pmb{z},t) \end{array} \right] \,,
\end{equation}
where the exponent $r=pq+1$ ensures the formal skew symmetry of the matrix-valued differential operator $\mathcal{J}$ (that is skew symmetry assuming zero boundary conditions for the arguments).
According to \eqref{eq:VDdef}, the time derivative of the energy functional (power balance) reads:
\begin{equation}
	\label{eq:PBE1}
	\dot{H}=\int_\Omega{ \delta_{{q}}H \wedge \dot{\alpha}^q + \delta_{{p}}H \wedge \dot{\alpha}^p}=\left< \delta_{{q}}H \right. \left| \dot{\alpha}^q \right>_\Omega + \left< \delta_{{p}}H \right. \left| \dot{\alpha}^p \right>_\Omega \,.
\end{equation}

According the state equations \eqref{eq:S2CL} and to the integration by part formula \eqref{eq:nppipp}, this power balance may be written:
\begin{equation}
	\label{eq:PBE2}
	\dot{H}= \left<  \delta_{{q}}H \right. \left| (-1)^p \delta_{{p}}H \right>_{\partial \Omega} = (-1)^{n-q} \int_{\partial \Omega}{  \left. \delta_{{q}}H \right|_{\partial \Omega}  \wedge \left. \delta_{{p}}H \right|_{\partial \Omega}} \,.
\end{equation}
This latter formula suggests of the following boundary port variables:
\begin{equation}
	\label{eq:BPVdef}
	\left[ \begin{array}{c} {e_\partial} \\  {f_\partial} \end{array} \right]  = \left[ \begin{array}{cc} (-1)^{n-q} \mbox{tr} & 0 \\ 0 & \mbox{tr} \end{array} \right] \left[ \begin{array}{c}  \delta_{{q}}H (\pmb{z},t) \\ \delta_{{p}}H (\pmb{z},t) \end{array} \right] \,,
\end{equation}
defined in such a way that the power balance equation \eqref{eq:PBE1} and \eqref{eq:PBE2} may be written:
\begin{equation}
	\label{eq:PBE3}
	\left< {e_q} \right. \left| {f^q} \right>_\Omega + \left< {e_p} \right. \left| {f^p} \right>_\Omega +
 \left<  {e_\partial} \right. \left| {f_\partial} \right>_{\partial \Omega} = 0 \,,
\end{equation}
where the flows variables ${f^q}\in L^2\Lambda^q(\Omega)$ and ${f^p}\in L^2\Lambda^p(\Omega)$ are defined as
\begin{equation}
	\label{eq:fqpdef}
	\begin{array}{c} {f^q}(\pmb{z},t):={-{\dot{\alpha}^q}(\pmb{z},t)} \,, \\   {f^p}(\pmb{z},t):={-{\dot{\alpha}^p}(\pmb{z},t)} \,.\end{array}
\end{equation}

The power balance equation \eqref{eq:PBE3} simply states that the time derivative of the energy increase inside the domain $\Omega$ equals the power supplied through the boundary $\partial \Omega$. As it may be seen, we have extended the systems of two conservation laws \eqref{eq:S2CL} with a boundary power supply and related boundary port variables, obtaining an open system of two conservation laws with boundary energy flows. The {\em explicit definition} of an open port-Hamiltonian system of two canonically interacting conservation laws is given by the distributed state equation \eqref{eq:S2CL} together with the definition of the boundary port variables \eqref{eq:BPVdef}. It leads to the structural power balance equation \eqref{eq:PBE3} which is independent of the specific considered Hamiltonian function (i.e. from the effort constitutive equations \eqref{eq:VDdef}). Many 1D, 2D or 3D examples, either linear or nonlinear may be recast in that framework and satisfy this definition \cite{VanderSchaft2014,duindam2009modeling}. Even parabolic systems may be formally represented  with the skew-symmetric operator $\mathcal{J}$ with a appropriate definition of the effort variables~\cite{Vu2016}.

\subsection{The irrotational 2D Shallow Water Equation example}
\label{subsec:wave2d}
We will consider as a running example for this paper the two-dimensional irrotational Shallow Water Equations (2D SWE) which describe the flow of an inviscid liquid where the horizontal components of the velocity field may be averaged on the water level and where the vertical velocity component may be omitted (low depth or shallow water assumption). Besides, we will consider a ``non rotating" flow. It is known that the corresponding 2D SWE express then the mass and momentum balance equations. Therefore, we will choose for the energy state variables the mass density (which is a 2-form proportional to the water level $h(\pmb{z},t)$) and the momentum density (which is a 1-form). For instance, using Cartesian coordinates, one would choose $\alpha^q:=h(\pmb{z},t) dx\wedge dy$ (where $h(\pmb{z},t)$ denotes the water level) and $\alpha^p:=\rho \left( u(\pmb{z},t) dx+ v(\pmb{z},t) dy\right) $ where $u(\pmb{z},t)$ and $v(\pmb{z},t)$ denote the horizontal components of the fluid velocity while $\rho$ denotes the fluid mass density. According to the previous notations, $\alpha^q\in L^2\Lambda^q(\Omega)$ is a 2-form and $\alpha^p\in L^2\Lambda^1(\Omega)$ a 1-form, both defined in the 2-dimensional ($n=2$) horizontal  spatial domain $\Omega$ of the flow. Using these energy state variables, one gets for the total (kinetic and potential) energy\footnote{Since the Hamiltonian density and the resulting constitutive equations relate differential forms of different degrees and are metric dependent, the Hodge star is introduced. We will consider in this paper the usual uniform metric.}  inside the domain $\Omega$:
\begin{equation}
	\label{eq:HSW2Ddef}
	H(\alpha^q,\alpha^p):= \int_\Omega{ \frac{\rho g\,(\star{\alpha^q}) \alpha^q}{2} + \frac{\star{\alpha^q}\left( \alpha^p \wedge \star{\alpha^p} \right)}{2 \rho} } \,.
\end{equation}
Therefore, the co-energy variables are defined as:
\begin{equation}
	\label{eq:eqepSW2Ddef}
	\begin{split}
		e_q  & = \delta_q H = \rho g (\star{\alpha^q}) + 
\frac{1}{2\rho} \star \left( (\star{\alpha^p}) \wedge \alpha^p \right) \,,\\
		e_p  & = \delta_p H = - \frac{(\star{\alpha^q}) 
(\star{\alpha^p})}{\rho}\,,
	\end{split}
\end{equation}
which are respectively the hydrodynamic pressure $e_q\in H^1\Lambda^{0}(\Omega)$ (a 0-form or intensive variable in the 2D domain $\Omega$) and the volume flow $e_p \in H^1\Lambda^{1}(\Omega)$ (which is indeed a 1-form in the 2D spatial domain). For instance, using the same Cartesian coordinates as previously, one gets $e_q= \rho g h + \frac{\rho}{2} \left(u^2+v^2\right)$ and $e_p= h \left( v\,dx -  u\,dy \right)$. Using these co-energy variables, the Hamiltonian system of two canonically interacting conservation laws \eqref{eq:S2CL} reads:
\begin{equation}
	\label{eq:2DSWECovForm}
	\left[ \begin{array}{c} \dot{\alpha}^q \\  \dot{\alpha}^p \end{array} \right]  = \left[ \begin{array}{cc} 0 & \mbox{d} \\ -\mbox{d} & 0 \end{array} \right] \left[ \begin{array}{c}  e_q \\ e_p \end{array} \right] \,,
\end{equation}
which are exactly the usual irrotational 2D SWE. For instance \eqref{eq:2DSWECovForm} reads in Cartesian coordinates using the usual vector calculus notations:
\begin{equation}
	\label{eq:2DSWECartCoord}
	\left[ \begin{array}{c} \dot{h} \\  \rho \left[ \begin{array}{c} \dot{u} \\ \dot{v} \end{array} \right] \end{array} \right]  = \left[ \begin{array}{cc} 0 & -\mbox{div} \\ -\mbox{grad} & 0 \end{array} \right] \left[ \begin{array}{c} \rho g h + \rho\frac{u^2+v^2}{2} \\ h \left[ \begin{array}{c} {u} \\ {v} \end{array} \right] \end{array} \right] \,.
\end{equation}

\subsection{The geometric structure of port-Hamiltonian systems}
\label{subsec:SDS}
In the previous section, we proposed a port-Hamiltonian formulation for open systems of two canonically interacting conservation laws (distributed state space equations \eqref{eq:S2CL} and the boundary equations \eqref{eq:BPVdef}). Since boundary energy flows are considered, boundary port variables are needed to derive the power balance equation \eqref{eq:PBE3}. The port-Hamiltonian model (\ref{eq:S2CL}, \ref{eq:BPVdef}) may then be implicitly defined as a linear subspace in the Bond space of effort and flow variables which embeds boundary effort and flow variables. In turn, this linear subspace may be geometrically defined as the linear subspace which is maximally isotropic with respect to some inner product associated to the natural power product - or power form in the Bond space - between effort and flow variables. Therefore, in the sequel, we aim at structure-preserving discretization which will preserve this power form in the approximation spaces. We will speak about symplectic discretization in the sense that this power form is preserved. In this section, we will define the Bond space, the power symplectic form, the associated inner product and the associated Stokes-Dirac structure which implicitely defines the port Hamiltonian model (\ref{eq:S2CL}, \ref{eq:BPVdef}). Readers are referred to \cite{VanderSchaft2002} for details about this representation.

Let the Bond space of extended flow and effort variables be $\mathcal{B}:=\mathcal{F}\times \mathcal{E}$ with
\begin{equation}
\label{ExFbondspacedef}
	\begin{split}
		\mathcal{F} & := L^2\Lambda^q(\Omega)\times  L^2\Lambda^p(\Omega)\times  L^2\Lambda^{n-q}(\partial \Omega) \\
		\mathcal{E} & := H^1\Lambda^{n-q}(\Omega)\times  H^1\Lambda^{n-p}(\Omega)\times  L^2\Lambda^{n-p}(\partial \Omega)
	\end{split}
\end{equation}

We may define on this Bond space the real power pairing or power form which maps any effort-flow vector $ \left( {e} , {f} \right) \equiv \left( (f^q,f^p,f_\partial),(e_q,e_p,e_\partial) \right) \in \mathcal{B}$ to
\begin{equation}
\label{powerformdef}
\left< {e} , {f} \right> := \left< {e_q} \right. \left| {f^q} \right>_\Omega + \left< {e_p} \right. \left| {f^p} \right>_\Omega +
 \left<  {e_\partial} \right. \left| {f_\partial} \right>_{\partial \Omega}
\end{equation}
in such a way that every pair $ \left( {e} , {f} \right)$ of extended effort and flow variables in the Bond space, satisfying the port Hamiltonian equations \eqref{eq:S2CL} and  \eqref{eq:BPVdef}), also satisfies the power  balance equation $ \left< {e} | {f} \right> =0$. From the power pairing \eqref{powerformdef}, we may define the following symmetric bilinear form:
\begin{equation}
\label{innerproductdef}
	\begin{split}
		\ll \cdot , \cdot \gg & :\mathcal{B}\times \mathcal{B} \rightarrow   \R \\
		& \ll  \left( {e_1} , {f_1} \right) , \left( {e_2} , {f_2} \right) \gg := \frac{1}{2} \left( \left< \left. {e_1} \right| {f_2} \right> + \left< \left. {e_2} \right| {f_1} \right> \right)
	\end{split}
\end{equation}
With the help of this symmetric bilinear form (inner product on $\mathcal{B}$) we may define the {\em Dirac structure associated to the power pairing} \eqref{powerformdef} as the linear subspace $\mathcal{D} \subset \mathcal{B}$ which is maximally isotropic, that is such that $\mathcal{D} = \mathcal{D}^\perp$ where the orthogonality is defined with respect to the inner product $\ll \cdot , \cdot \gg$. In particular, any $ \left( {e} , {f} \right) \in \mathcal{D}$ satisfies $\ll  \left( {e} , {f} \right) , \left( {e} , {f} \right) \gg = 0$, hence the power balance $\left< {e} , {f} \right> =0$. 

Dirac interconnection structure may used to define implicitely the dynamics of port Hamiltonian systems. In particular, the linear subspace $\mathcal{D}$ in the Bond space $\mathcal{B}:=\mathcal{F}\times \mathcal{E}$, with $\mathcal{F}$ and $\mathcal{E}$ as in \eqref{ExFbondspacedef}, which is defined by:
\begin{equation}
\label{implPHSdef}
\mathcal{D}:=\left\{  
\left( ({f^q},{f^p},f_\partial), ({e_q},{e_p}, e_\partial) \right)\in \mathcal{B} 
\left| 
\begin{array}{l}
\left[ \begin{array}{c} f^q \\  f^p \end{array} \right]  = \left[ \begin{array}{cc} 0 & - \mbox{d} \\ \mbox{d} & 0 \end{array} \right] 
\left[ \begin{array}{c}  e_q \\ e_p \end{array} \right] 
\\
\left[ \begin{array}{c} {e_\partial} \\  {f_\partial} \end{array} \right]  = \left[ \begin{array}{cc} (-1)^{p+1} \mbox{tr} & 0 \\ 0 & \mbox{tr} \end{array} \right] \left[ \begin{array}{c}  e_q \\ e_p \end{array} \right]
\end{array} 
\right.
\right\}
\end{equation}
is a Dirac structure associated to the natural power pairing \eqref{powerformdef}. This is proved by using the generalized Stokes theorem \cite{VanderSchaft2002}. Therefore, in this particular case, the interconnection structure is called a Stokes-Dirac structure. The dynamics \eqref{eq:S2CL}, generated by the Hamiltonian function $H(\alpha^q,\alpha^p)$ (see definition \eqref{eq:Hdef}),  with boundary energy flow and  port boundary variables \eqref{eq:BPVdef}, may be implicitly defined by:
\begin{equation}
\label{impldyndef}
\left( (-\dot{\alpha}^q,-\dot{\alpha}^p,f_\partial), (\delta_qH,\delta_pH, e_\partial) \right) \in \mathcal{D}
\end{equation}
In that sense, we will say that \eqref{impldyndef} defines a boundary port Hamiltonian system of two canonically interconnected conservation laws. In order to define a well-posed Cauchy problem, boundary conditions still need be chosen either for $e_\partial$ or $f_\partial$.  
\begin{remark}
The Stokes-Dirac structure \eqref{implPHSdef}, associated to the natural power bilinear form \eqref{powerformdef},
may be used to represent hyperbolic systems, either linear (when the Hamiltonian density is quadratic) or nonlinear (in the other cases). It may even be used to represent parabolic systems when the effort differential forms do not derive from the same Hamiltonian operator \cite{baaiu2009port}. The choice of Dirichlet boundary conditions for $e_\partial$ or $f_\partial$ will lead to a well-posed system. They are however many other possible choices for admissible boundary conditions which lead to well-posed systems. For linear systems (quadratic Hamiltonian), these admissible boundary conditions have been parameterized elegantly \cite{Gorrec2005}. Many examples of physical systems have been represented using the port Hamiltonian formulation and its implicit representation using Stokes-Dirac structures, including Maxwell field equations and Navier-Stokes flow problems \cite{VanderSchaft2002},  beam and membrane equations \cite{duindam2009modeling}, vibro-acoustic problems \cite{trenchant2015port}, shallow water flow problems
\cite{hamroun2006port,pasumarthy2008port,Cardoso-Ribeiro2015a}, advection-diffusion or adsorption problems \cite{baaiu2009port}, tokamak plasma MHD problems \cite{Vu2016},etc.
\end{remark}

\subsection{The partitioned weak form for port-Hamiltonian systems of two conservation laws}
\label{EDFWF}

We will follow the approach presented in section \ref{subsec:wfSWE} for the 1D SWE and generalize it to the general covariant formulation \eqref{eq:S2CL} for port Hamiltonian systems of two conservation laws. Let $v_q \in H^1 \Lambda^{n-q} \left( \Omega \right)$  and $v_p \in L^2\Lambda^{n-p} \left( \Omega \right)$ denote any arbitrary test differential forms. We may obtain from the strong formulation \eqref{eq:S2CL} the following weak form:
\begin{equation}
	\begin{split}
		-\left< v_q,\dot{\alpha}^q \right>_\Omega & = - \left< v_q , \mathrm{d} e_p \right>_\Omega \\ 		
		-\left< v_p,\dot{\alpha}^p \right>_\Omega & = (-1)^{pq} \left< v_p , \mathrm{d} e_q \right>_\Omega
	\end{split}
 \end{equation}
Integrating by part the first state equation only, and using the integration by part formula \eqref{eq:nppipp}, we get the following partitioned weak form:
\begin{equation}
\label{eq:partweakcovform}
	\begin{split}
		-\left< v_q,\dot{\alpha}^q \right>_\Omega & = (-1)^{p} \left(  \left< \mathrm{d}v_q , e_p \right>_\Omega - \left< v_q , e_p \right>_{\partial \Omega} \right) \\ 		
		-\left< v_p,\dot{\alpha}^p \right>_\Omega & = (-1)^{pq}  \left< v_p , \mathrm{d} e_q \right>_\Omega
	\end{split}
 \end{equation}

\begin{remark}
With the particular choice $v_q =e_q$ and $v_p = e_p$, one gets:
\begin{equation}
\label{eq:partweakcovorm_powerbalance}
	\begin{split}
		-\left< e_q,\dot{\alpha}^q \right>_\Omega & = (-1)^{p} \left(  \left< \mathrm{d}e_q , e_p \right>_\Omega - \left< e
_q , e_p \right>_{\partial \Omega} \right) \\ 		
		-\left< e_p,\dot{\alpha}^p \right>_\Omega & = (-1)^{pq}  \left< e_p , \mathrm{d} e_q \right>_\Omega
	\end{split}
 \end{equation}
Therefore the power-balance equation \eqref{eq:Hdot1D} reads:
\begin{equation}
		\dot{H}   = \left< e_q,\dot{\alpha}^q \right>_\Omega + \left< e_p,\dot{\alpha}^p \right>_\Omega =  (-1)^p \left< e_q , e_p \right>_{\partial \Omega}
\end{equation}
which shows that the power balance is also preserved in the weak formulation.
\end{remark}
\begin{remark}
	 Instead of integrating the first equation by parts in  \eqref{eq:partweakcovform}, we could integrate the second equation, which would also lead to another skew-symmetric structure; the boundary inputs and outputs would not be the same either.
\end{remark}
In the case of the irrotational 2D SWE, written in Cartesian coordinates, using usual vector calculus notations, the partitioned weak formulation \eqref{eq:partweakcovform} transforms into
\begin{equation}
	\label{eq:weakformpartitioned2D}
	\begin{split}
		\int_\Omega {v_q {\dot{\alpha}_q}} \dx \dy & = \int_\Omega { \left( \pmb{\nabla}v_q \right) \cdot \pmb{e}_{\pmb{p}}}\dx\dy - \int_{\partial\Omega} v_q \pmb{n} \cdot \pmb{e}_{\pmb{p}} \ds\\
		\int_\Omega {\pmb{v_p} \cdot \pmb{\dot{\alpha}_p}}\dx\dy & = - \int_\Omega{\pmb{v_p} \cdot  \pmb{\nabla} e_q}\dx\dy \,,
	\end{split}
\end{equation}
In this latter equation, notations have been introduced both for simplification and to distinguish between scalar coordinates (for functions and 2-forms) and vector coordinates (for 1-forms). For instance, ${\alpha}_q$ has been used instead of $\star{\alpha}_q=h(t,x,y)$ and denotes the water level (proportional to the water surfacic density), while $v_q$ and $e_q$ directly stands for the corresponding functions (0-forms). The bold vector notations $\pmb{{\alpha}_p}:=\rho h(t,x,y)[u(t,x,y),v(t,x,y)]^T$, $\pmb{e}_{\pmb{p}}$ and $\pmb{v}_{\pmb{p}}$ denote respectively the vectors of Cartesian coordinates for the momentum (${\alpha}^p$), water flow ($e_p$) and test ($v_p$) which are all 1-forms. The notation $\pmb{\nabla}$ denotes the usual 2D Nabla differential operator while $\pmb{n} \ds$ is the vector of Cartesian coordinates for the length 1-form (perpendicular to the boundary $\partial \Omega$). Note that the term $- \pmb{n} \cdot  \pmb{e}_{\pmb p}$ denotes the boundary port variable $f_\partial$ (perpendicular water flow at the boundary) which will be chosen as the boundary input $u_\partial$.

In the following section, we discretize the partioned weak-form \eqref{eq:weakformpartitioned2D}, and we show that the resulting system is a finite-dimensional port-Hamiltonian system preserving the power-balance of the original system.
}

\section{The Partitioned Finite-Element Method (PFEM) for the 2D SWE}
\label{sec:2DPFEM}
In this section, the partitioned weak-form representation for the Shallow Water Equations, presented in \S~\ref{EDFWF}, will be projected into finite-dimensional approximation spaces in such a way as to preserve the total mass and momentum conservation equations and the underlying Dirac structure of the original port-Hamiltonian model. This section is divided in two parts, firstly, the weak-form is discretized in \S~\ref{subsec:discrete2D} and a finite-dimensional port-Hamiltonian system is obtained. Secondly, we show how to obtain the discrete constitutive relationships in \S~\ref{subsec:nonlinearconstitutive}. 

\subsection{Structure-preserving finite element discretization}
\label{subsec:discrete2D}

Let us approximate the energy variables  $\alpha_q(x,y,t)$ using the following basis with $N_q$ elements:
\begin{equation}
	\label{eq:approxvariables2D_q}
	\begin{split}
	\alpha_q(x,y,t) \approx \alpha_q^{ap}(x,y,t) & := \sum_{i=1}^{N_q} \phi^i_q(x,y) \alpha_q^i(t)  = \pmb{\phi}_q(x,y)^T \pmb{\alpha}_q(t) \,,
	\end{split}
\end{equation}
where $\phi_q^i(x,y)$, $i \in \{1,\dots, N_q\}$ are the chosen approximation basis functions in $H^1(\Omega)$, and $\alpha^i_q(t)$ are the approximation coordinates for $\alpha_q^{ap}(x,y,t)$. The test functions $v_q$  and the co-energy variables $e_q$ are also approximated using $\pmb{\phi}_q(x,y)$.

Similarly, the vectorial variable $\pmb{\alpha_p}(x,y,t)$ is approximated as:
\begin{equation}
	\label{eq:approxvariables2D_p}
	\begin{split}
	\pmb{\alpha_p}(x,y,t) \approx \pmb{\alpha_p}^{ap}(x,y,t) & := \sum_{i=1}^{N_p} \pmb{\phi_p}^i(x,y) \alpha_p^i(t) = \pmb{\Phi}_p(x,y)^T \pmb{\alpha}_p(t)\,,
	\end{split}
\end{equation}
where ${\pmb{\phi_p}^i(x,y) = \left[\begin{matrix} \phi_p^{x,i} (x,y) \\ \phi_p^{y,i} (x,y) \end{matrix}\right]}$ represents a 2D-vectorial basis function and, consequently, $\pmb{\phi_p}(x,y)$ is an $N_p \times 2$ matrix. The variables $\alpha_p^i(t)$ are the approximation coordinates for $\alpha_p^{ap}(x,y,t)$ in the $\pmb{\Phi}_p(x,y)$ approximation space. The test functions $\pmb{v}_p$ and the co-energy variable $\pmb{e}_p$ are also approximated using $\pmb{\Phi}_p(x,y)$.

\begin{remark}
  \label{rem6}
 Note that a particular choice for $\pmb{\Phi}_p$ is :
\begin{equation}
	\pmb{\Phi}_p(x,y) = \left[\begin{matrix} \pmb{\phi}_{p}^x(x,y) & 0 \\ 0 & \pmb{\phi}_{p}^y(x,y) \end{matrix} \right]\,
\end{equation}
such that we can decompose the variables with index $p$ in their Cartesian components as: $\pmb{\alpha_p} = \left[\begin{matrix} \pmb{\alpha}_{p}^x \\ \pmb{\alpha}_{p}^y\end{matrix} \right]$.

	With this particular choice, we recover the case where the basis functions components of the vectorial variables are decoupled. We studied this case in \cite{CardosoRibeiro2018}.
 \end{remark}
 
 Finally, the boundary input can be discretized using any one-dimensional set of basis functions, say $\pmb{\psi}(s)=[\psi^i(s)]$:
 \begin{equation}
	\label{eq:boundary2ddiscrete}
	u_\partial(s,t) \approx u^{ap}_\partial(s,t) := \sum_{i=1}^{N_\partial} \psi^i(s) u_{\partial}^i(t) = \pmb{\psi}(s)^T \pmb{u}_{\partial}(t)\,.
 \end{equation} 

 \begin{remark}
	In the sequel, in the implementation of our finite element method, we conveniently chose $\pmb{\psi}(s)$ as $\pmb{\phi}_q(x(s),y(s))$ evaluated on the boundary. Other choices could be investigated.
\end{remark}
 
 From the substitution of the approximated variables, \eqref{eq:approxvariables2D_q},\eqref{eq:approxvariables2D_p} and \eqref{eq:boundary2ddiscrete}, in the weak form \eqref{eq:weakformpartitioned2D}, the following finite-dimensional equations are obtained:
 \begin{equation}
	\label{eq:finidemensionaleqs0}
	\begin{split}
	\pmb{v_q}^T \left[\int_\Omega{\pmb{\phi}_q \pmb{\phi}_q^T \dx\dy} \right]\, \dot{\pmb{\alpha}}_q  = & \pmb{v_q}^T \left[\int_\Omega{ \left[ \begin{matrix} \frac{\partial \pmb{\phi}_{q}}{\partial x} & \frac{\partial \pmb{\phi}_{q}}{\partial y} \end{matrix} \right]\pmb{\Phi}_p^T \dx\dy} \right] \, \pmb{e_p} + \\
		& - \pmb{v_q}^T \left[\int_{\partial \Omega}{\pmb{\phi}_q(x(s),y(s)) \Psi^T(s) \ds  }\right]\, \pmb{u}_\partial(t) \,,\\
	\pmb{v_p}^T\left[\int_\Omega{\pmb{\Phi}_p \pmb{\Phi}_p^T \dx\dy}\right]\, \dot{\pmb{\alpha}}_p = &  - \pmb{v_p}^T \left[\int_\Omega{\pmb{\Phi}_p \left[ \begin{matrix} \frac{\partial\pmb{\phi}_{q}^T}{\partial x} \\ \frac{\partial\pmb{\phi}_{q}^T}{\partial y} \end{matrix} \right] \dx\dy}\right]\, \pmb{e}_q \,.
	\end{split}
 \end{equation}
Since these equations should remain valid for any test functions coordinates $\pmb{v_q}$ and $\pmb{v_p}$, one gets:
\begin{equation}
	\label{eq:finitedimensionaleqs}
	\begin{split}
		{M_q} \dot{\pmb{\alpha}}_q (t)  = & {D} \pmb{e}_p (t)+ {B} \pmb{u}_\partial(t) \,,\\
		{M_p} \dot{\pmb{\alpha}}_p (t)= &  - {D^T} \pmb{e}_q (t)\,,
	\end{split}
 \end{equation}
 where $M_q$ and $M_p$ are square mass matrices (of size $N_q \times N_q$ and $N_p \times N_p$, respectively), defined as:
 \begin{equation}
	M_q := \int_\Omega{\pmb{\phi}_q \pmb{\phi}_q^T \dx\dy}\,,\quad\quad\quad
	M_p := \int_\Omega{\pmb{\Phi}_p \pmb{\Phi}_p^T \dx\dy}\,.
 \end{equation}
	 The matrix $D$ is of size $N_q \times N_p$, defined as:
	 \begin{equation}
		D := \int_\Omega{ \left[ \begin{matrix} \frac{\partial \pmb{\phi}_{q}}{\partial x} & \frac{\partial \pmb{\phi}_{q}}{\partial y} \end{matrix} \right]\pmb{\Phi}_p^T \dx\dy}\,.
	 \end{equation}
	 and $B$ is an $N_q \times N_\partial$ matrix:
\begin{equation}
	B := \int_{\partial \Omega}{\pmb{\phi}_q(x(s),y(s)) \pmb{\psi}^T(s) \ds  }\,.
\end{equation}

 Defining $\pmb{y}_\partial(t)$, the output conjugated to the input $\pmb{u}_\partial(t)$ as:
 \begin{equation}
	\pmb{y}_\partial(t) := B^T \pmb{e}_q(t)\,,
 \end{equation}
 the approximated system can be written using the following finite-dimensional Dirac structure representation:
 \begin{equation}
	\label{eq:finitedim_2dwaveeq}
	\begin{split}
		\left[ \begin{matrix} M_q & 0 \\ 0 & M_p \end{matrix}\right] \left[ \begin{matrix} {\pmb{f}}_q(t) \\ {\pmb{f}}_p(t) \end{matrix}\right] & = \left[ \begin{matrix} 0 & -D \\ D^T & 0 \end{matrix}\right] \left[ \begin{matrix} {\pmb{e}}_q(t) \\ {\pmb{e}}_p(t) \end{matrix}\right] + \left[ \begin{matrix} -B \\ 0 \end{matrix}\right] \pmb{u}_\partial(t) \\
		\pmb{y}_\partial(t) & = \left[ \begin{matrix} B^T & 0 \end{matrix}\right] \left[ \begin{matrix} {\pmb{e}}_q(t) \\ {\pmb{e}}_p(t) \end{matrix}\right]
	\end{split}
\end{equation}
where ${\pmb{f}}_q(t) := -\dot{\pmb{\alpha}}_q(t)$ and ${\pmb{f}}_p(t) := -\dot{\pmb{\alpha}}_p(t)$ denote the vector coordinates for the flow approximations coordinates in the $\phi^i_q(x,y)$  and $\pmb{\phi_p}^i(x,y)$ approximation bases, that is $f_q^{ap}(x,y,t) = \pmb{\phi}_q(x,y)^T \pmb{f}_q(t)$ and $f_p^{ap}(x,y,t) = \pmb{\Phi}_p(x,y)^T \pmb{f}_p(t)$.

 From the definition of the co-energy variables as the variational derivatives of the Hamiltonian with respect to $\alpha_q(x,y,t)$ and $\pmb{\alpha}_p(x,y,t)$, the time-derivative of the continous Hamiltonian is given by:
 \begin{equation}
	\begin{split}
		\dot{H} & = \int_{\Omega}{ \left(\dot{\pmb{\alpha}}_{\pmb{p}}(\pmb{x},t) \cdot \pmb{e}_{\pmb{p}}(\pmb{x},t)  +  \dot{\alpha}_q(\pmb{x},t) e_q(\pmb{x},t) \right)\dOmega }\,,
	\end{split}
\end{equation}
Using the approximations for the energy and co-energy variables, this power balance can be approximated as:
	\begin{equation}
	   \label{eq:Hdot2d_discreteCho}
	   \begin{split}
		\dot{H}_d   &= \int_{\Omega}{\left(\dot{\pmb{\alpha}}_p(t)^T \pmb{\Phi}_p(x,y)  \pmb{\Phi}_p(x,y)^T \pmb{e}_p(t)  + \dot{\pmb{\alpha}}_q(t)^T \pmb{\phi}_q(x,y) \pmb{\phi}_q(x,y)^T \pmb{e}_q(t)\right) \dOmega} \,,\\
		 & =  \dot{\pmb{\alpha}}_p (t)^T M_p \pmb{e}_p(t)  + \dot{\pmb{\alpha}}_q(t)^T M_q \pmb{e}_q(t) \,.
	   \end{split}
	\end{equation}

	The approximated equations \eqref{eq:finitedim_2dwaveeq}, together with the power balance \eqref{eq:Hdot2d_discreteCho}, provides a Dirac structure representation that is a projection of the Dirac-Stokes structure. As we did for the 1D SWE in \S~\ref{subsec:wfSWE}, we will now derive the corresponding  approximation of the port-Hamiltonian dynamics by restricting the Hamiltonian functional to the same approximation spaces. New energy and co-energy variables will be defined, so that we can find an explicit port-Hamiltonian representation for the approximated system.

	Since $M_q$ and $M_p$ are symmetric positive definite matrices, we can use Cholesky decomposition\footnote{The use of Cholesky is advantageous since it simplifies the numerical inversion of the mass matrices. Note that the procedure presented in this section slightly differs from the one presented in \S~\ref{subsec:wfSWE}, since different energy and co-energy variables are defined here.}. There exist triangular matrices $L_q$ and $L_p$, such that: $M_q = L_q L_q^T$ and $M_p = L_p L_p^T$. The power balance becomes:
	\begin{equation}
		\label{eq:powerbalancediscretized2d}
		\begin{split}
		 \dot{H}_d   & =  \left(L_p^T\,\dot{\pmb{\alpha}}_p (t)\right)^T  L_p \pmb{e}_p(t)  + \left(L_q^T \,\dot{\pmb{\alpha}}_q(t)\right)^T L_q \pmb{e}_q(t) \,,\\
		 & =  \dot{\tilde{\pmb{\alpha}}}_p (t)^T  \pmb{\tilde{e}}_p(t)  + \dot{\pmb{\tilde{\alpha}}}\,\,_q(t)^T \pmb{\tilde{e}}_q(t) \,.
		\end{split}
	 \end{equation}
	 where new energy and co-energy variables:
	 \begin{equation}
		\label{eq:newenergyvariables2dCho}
		\begin{split}
		   {\tilde{\pmb{\alpha}}}_p (t):=& L_p^T {\pmb{\alpha}}_p(t)\,, \quad  \quad \quad
		   {\tilde{\pmb{\alpha}}}_q (t):= L_q^T {\pmb{\alpha}}_q(t)\,, \\
		   {\tilde{\pmb{e}}}_p (t):=& L_p^T {\pmb{e}}_p(t)\,, \quad  \quad \quad
		   {\tilde{\pmb{e}}}_q (t):= L_q^T {\pmb{e}}_q(t)\,, \\
		\end{split}
	\end{equation}
   have been defined. Therefore, in order to write the power-balance as the total time derivative of the discrete Hamiltonian written as a function of the finite-dimensional vector coordinates for the energy variables, the newly defined reduced co-energy variables must be equal to the gradient of the Hamiltonian with respect to the energy variables:
   \begin{equation}
	   \begin{split}
		   \pmb{\tilde{e}}_q = \frac{\partial H_d}{\partial{\tilde{\pmb{\alpha}}}_q} \,, \\
		   \pmb{\tilde{e}}_p = \frac{\partial H_d}{\partial{\tilde{\pmb{\alpha}}}_p} \,,
	   \end{split}
   \end{equation}
   where the approximated Hamiltonian is defined as:
   \begin{equation}
   \label{eq:discretizedHamiltonianCho}
	   \begin{split}
		   H_d(\tilde{\pmb{\alpha}}_q,\tilde{\pmb{\alpha}}_p) := H&\left[\alpha_q(\pmb{x},t)= L_q^{-T} \tilde{\pmb{\alpha}}^T_q(t) \pmb{\phi_q}(\pmb{x}) \,, \right. \\
		   &\left.\alpha_p(\pmb{x},t)= L_p^{-T} \tilde{\pmb{\alpha}}_p^T(t) \pmb{\Phi_p}(\pmb{x})  \right]\,.
	   \end{split}
   \end{equation}
   
   Rewriting the finite-dimensional equations \eqref{eq:finitedimensionaleqs}, we get the following explicit finite-dimensional port-Hamiltonian system:
   \begin{equation}
   \label{eq:fin_dim_2d_phs_matrixCho}
	   \begin{split}
		   \left[
			   \begin{matrix}
			   \dot{\tilde{\pmb{\alpha}}}_q \\ \dot{\tilde{\pmb{\alpha}}}_p
			   \end{matrix}
		   \right]
		   =&
		   {\left[
			   \begin{matrix}
			   0 & L_q^{-1} D L_p^{-T} \\- L_p^{-1} D^T L_q^{-T} &0 
			   \end{matrix}
		   \right]}
		   \left[
			   \begin{matrix}
				   \pmb{\tilde{e}}_q \\
				   \pmb{\tilde{e}}_p 
			   \end{matrix}
		   \right]
		   + 
		   {\left[
			   \begin{matrix}
				L_q^{-1}B \\ 0
			   \end{matrix}
		   \right]} \pmb{u}_\partial \,,\\
		   \pmb{y}_\partial =& B^T L_q^{-T} \pmb{\tilde{e}}_q \,,
	   \end{split}
	\end{equation}
   where $\pmb{y}_\partial$ is the conjugated output of the discretized system.
   
   From \eqref{eq:Hdot2d_discreteCho} and \eqref{eq:fin_dim_2d_phs_matrixCho}, the time derivative of the approximated Hamiltonian is given by:
   \begin{equation}
	 \label{eq:discretepowerbalanceCho}
	   \begin{split}
		   \dot{H}_d(t) & = \dot{\tilde{\pmb{\alpha}}}_q^T \, \pmb{\tilde{e}}_q + \dot{\tilde{\pmb{\alpha}}}_p^T \, \pmb{\tilde{e}}_p \,,\\
			& = \left(\pmb{e}_p^T  L_p^{-1} D^T L_q^{-T} + \pmb{u}_\partial^T L_q^{-1}B \right) \, \pmb{\tilde{e}}_q -\pmb{\tilde{e}}_q^T L_q^{-1} D L_p^{-T}\, \pmb{\tilde{e}}_p   \,,\\
					& =  \pmb{y}^T_\partial  \, \pmb{u}_\partial\,.
	   \end{split}
   \end{equation}
   Note that  $\pmb{u}_\partial^T \pmb{y}_\partial$ is the discrete analog of the continuous power-balance equation \eqref{eq:powerbalance2d}. Furthermore, this power-balance is exactly preserved in the finite-dimensional approximation spaces. From the definition of the $B$ matrix~\eqref{eq:finidemensionaleqs0}, the definition of the approximated boundary input $u^{ap}_\partial(s,t) := \pmb{\psi}(s)^T \pmb{u}_{\partial}(t)$ and approximated co-energy variable $e_p^{ap}(x(s),y(s),t) := \pmb{\phi}_q(x(s),y(s))^T \pmb{e}_{q}(t)$, we get:
   \begin{equation}
	   \begin{split}
				   \dot{H}_d  & = \pmb{e}^T_q \, \int_{\partial \Omega}{\pmb{\phi}_q(x(s),y(s)) \pmb{\psi}^T(s) \ds  }   \, \pmb{u}_\partial\,, \\
					 & =  \, \int_{\partial \Omega}{e_q^{ap}(x(s),y(s),t) u_\partial^{ap}(s,t) \ds  }   \, .
	   \end{split}
   \end{equation}

   \begin{remark}
     \label{rem8}
	Note that using classical finite-elements 1D discretization basis for the boundary input, the coordinates $\pmb{u}_\partial(t)$ provide the values of the inflow ($-\pmb{n} \cdot \pmb{e}_{\pmb{q}}$) at the boundary nodes. For instance, in the case of shallow water equations, these are the values of volumetric influx into the system.
	The conjugated output $\pmb{y}_\partial$ is related with the curve integral of $e_q(x(s),y(s),t)$ along the elements. The co-energy variable $e_q(x(s),y(s),t)$ is the pressure, thus the discretized outputs coordinates $\pmb{y}_\partial(t)$ are related to the forces per unit length applied along the external boundary.

	Note that we can also define an output that is given by the point-wise values of the co-energy variables $e_q(s,t)$ evaluated on the boundary. In this case, a convenient choice of basis function for the approximation would be the same as the input, i.e.:
	\begin{equation}
		y^{ap}(s,t) = \pmb{\psi}^T(s) \hat{\pmb{y}}_\partial (t) \,,
	\end{equation}
	where the coordinates $\hat{\pmb{y}}_\partial(t)$ represents the values of pressure at the boundary nodes. The power-balance through the boundary is computed as:
	\begin{equation}
		\label{eq:powerbalance-pointwiseY}
		\begin{split}
			\dot{H}_d &= \int_{\partial \Omega}{ y_\partial^{ap}(s,t) y_\partial^{ap}(s,t) \ds} \\
			&= \hat{\pmb{y}}_\partial^T (t)\left(\int_{\partial \Omega}{ \pmb{\psi}(s)  \pmb{\psi}^T(s) \ds} \right){\pmb{u}}_\partial \\
			& = \hat{\pmb{y}}_\partial^T (t) M_\psi {\pmb{u}}_\partial\,,
		\end{split}
	\end{equation}
	where
	\begin{equation}
		M_\psi = \int_{\partial \Omega}{ \pmb{\psi}(s) \pmb{\psi}^T(s) \ds} \,
	\end{equation}
	is a symmetric positive-definite $N_\partial \times N_\partial$ mass matrix.

	Furthermore, since the power-balances \eqref{eq:discretepowerbalanceCho} and \eqref{eq:powerbalance-pointwiseY} must coincide, the following relationship between these two output definitions must hold:
	\begin{equation}
		\hat{\pmb{y}}_\partial = M_\psi^{-1} {\pmb{y}}_\partial\,.
	\end{equation}

	Consequently, from \eqref{eq:powerbalancediscretized2d} and \eqref{eq:powerbalance-pointwiseY}, the following power product must hold (satisfying \eqref{eq:finitedim_2dwaveeq}):
	\begin{equation}
		\label{powerformdef_findim}
	\left< \left. \left[\begin{matrix}\pmb{f}_p\\ \pmb{f}_p \\ \pmb{\hat{y}}_\partial \end{matrix}\right] \right| \left[ \begin{matrix}\pmb{e}_q\\ \pmb{e}_p \\ \pmb{{u}}_\partial \end{matrix}\right] \right> := \pmb{f}_p^T M_p e_p  + \pmb{f}_q^T M_q \pmb{e}_q  +   \pmb{u}_\partial^T M_\psi \pmb{\hat{y}}_\partial = 0 \,.
	\end{equation}
	Alternatively, the power balance can be reduced to (satisfying \eqref{eq:fin_dim_2d_phs_matrixCho}):
	\begin{equation}
		\label{powerformdef_findim}
	\left< \left. \left[\begin{matrix}\pmb{\tilde{f}}_q\\ \pmb{\tilde{f}}_p \\ \pmb{{y}}_\partial \end{matrix}\right] \right| \left[ \begin{matrix}\pmb{\tilde{e}}_q\\ \pmb{\tilde{e}}_p \\ \pmb{{u}}_\partial \end{matrix}\right] \right> := \pmb{\tilde{f}}_p^T \pmb{\tilde{e}}_p  + \pmb{\tilde{f}}_q^T \pmb{\tilde{e}}_q  +   \pmb{u}_\partial^T \pmb{{y}}_\partial = 0 \,,
	\end{equation}
	where $\pmb{\tilde{f}}_q (t) := -\dot{\pmb{\tilde{\alpha}}}\,\,_q (t)$ and $\pmb{\tilde{f}}_p (t):= {-\dot{\pmb{\tilde{\alpha}}}\,\,}_p (t)$.

	With the help of a symmetric bilinear form as (\ref{innerproductdef}), using any of the previous power products, finite-dimensional Dirac interconnection structures can be defined.
	
\end{remark}

\subsection{Obtaining the non-linear constitutive relationships: discretization of the Hamiltonian}
\label{subsec:nonlinearconstitutive}
In the previous section, a finite-dimensional Dirac structure was obtained for the 2D Shallow Water Equations, relating the energy and co-energy variables as well as the boundary inputs and outputs. The next step is to obtain the constitutive relationships from the Hamiltonian. 

The Hamiltonian of the 2D SWE \eqref{eq:HSW2Ddef} can be rewritten using the coordinate variables as:
\begin{equation}
	H[\alpha_q(x,y,t), \pmb{\alpha_p}(x,y,t)] := \frac{1}{2}\int_\Omega{\left( \frac{{\alpha_q(x,y,t)} \|\pmb{\alpha_p}\|^2}{\rho} + \rho g({\alpha_q}(x,y,t))^2   \right) }\dOmega
\end{equation}

The energy variables are restricted to the approximation spaces for $\alpha_q(x,y,t)$ and $\pmb{\alpha_p}(x,y,t)$. From \eqref{eq:discretizedHamiltonianCho}, the discretized Hamiltonian reads:
\begin{equation}
\label{eq:2DSWEconstrel}
	{H}_d(\tilde{\pmb{\alpha}}_q,\tilde{\pmb{\alpha}}_p):=\frac{1}{2}\int_{\Omega} \left( \frac{\pmb{\phi}^T_q(x,y) L_q^{-1} \tilde{\pmb{\alpha}_q} \left( \pmb{\Phi}^T_p(x,y) L_p^{-1} \tilde{\pmb{\alpha}_p} \right)^2}{\rho} + \, {\rho g} \left( \pmb{\phi}^T_q(x,y) L_q^{-1} \tilde{\pmb{\alpha}_q} \right)^2 \right) \dOmega
\end{equation}
The reduced effort variables are obtained from the gradient of the discretized Hamiltonian:
\begin{equation}
	\label{eq:2DSWEconstrel4}
	\begin{split}
		\frac{\partial {H}_d}{\partial \tilde{\pmb{\alpha}}_q} & = {\rho g}L_q^{-T}\tilde{\pmb{\alpha}}_q(t)+\left[  \frac{1}{2\rho} \int_{\Omega}  L_q^{-T}  \pmb{\phi}_q(x,y) \tilde{\pmb{\alpha}}_p^T(t)  \pmb{\Phi}_p(x,y)\pmb{\Phi}^T_p(x,y) L_p^{-1}  \dOmega \right] \tilde{\pmb{\alpha}}_p(t)\\
		\frac{\partial {H}_d}{\partial \tilde{\pmb{\alpha}}_p} & = \tilde{\pmb{\alpha}}_q^T(t) \left[  \frac{1}{\rho}\int_{\Omega}  L_q^{-T}  \pmb{\phi}_q(x,y)  L_p^{-T} \pmb{\Phi}_p(z) \pmb{\Phi}_p^T(z)  L_p^{-1}   \dOmega \right] \tilde{\pmb{\alpha}}_p(t)
	\end{split}
\end{equation}
Both constitutive equations exhibits nonlinear terms. In order to compute them, the following procedured was used. The first equation can be written as:
\begin{equation}
	\frac{\partial \tilde{H}}{\partial \tilde{\pmb{\alpha}}_q} = L_q^{-T} \underbrace{\left({\rho g}\tilde{\pmb{\alpha}}_q(t)+\left[  \frac{1}{2\rho} \int_{\Omega}   \pmb{\phi}_q(x,y) \tilde{\pmb{\alpha}}_p^T(t)  \pmb{\Phi}_p(x,y)\pmb{\Phi}^T_p(x,y)  \dOmega \right] \tilde{\pmb{\alpha}}_p(t) \right)}_{\pmb{w}(t)}
\end{equation}
where the components of $\pmb{w}(t)$ can be computed as:
\begin{equation}
	w_i(t) = {\rho g} \tilde{{\alpha}}_{q,i}(t)+ \tilde{\pmb{\alpha}}_p^T(t) \left(\int_{\Omega}   {\phi}_{q,i}(x,y)   \pmb{\Phi}_p(x,y)\pmb{\Phi}^T_p(x,y)  \dOmega \right) \tilde{\pmb{\alpha}}_p(t) 
\end{equation}
Note that $\int_{\Omega}   {\phi}_{q,i}(x,y)   \pmb{\Phi}_p(x,y)\pmb{\Phi}^T_p(x,y)  \dOmega$, for $\{i = 1, \dots, N_q\}$  are $N_q$ matrices of dimension $N_p \times N_p$, that can be computed once and remains constant. Similarly, the second constitutive relationship can be written as a function of constant matrices.

\section{Numerical experiments}
\label{sec:numericalresults}

In this section, we present numerical experiments to test the Partitioned Finite Element Method. Firstly, results for the 1D SWE are presented in \S~\ref{ssNumRes1D}. Then, the 2D case is presented in \S~\ref{ssNumRes2D}.

\subsection{One-dimensional Shallow Water Equations}
\label{ssNumRes1D}

The example presented in Section \ref{sec:introexample} was implemented using finite elements with polynomial basis functions.

Firstly, a spectral convergence analysis of the numerical method was done. The eigenvalues obtained from the linearized numerical model were compared to the exact eigenfrequencies of the linear wave equation with constant coefficients. The inputs ${u}_\partial(t)$ of \eqref{eq:0DSWEPCH} were considered to be zero.

\begin{figure}[h]
	\begin{center}
	   \includegraphics[width=0.7\textwidth]{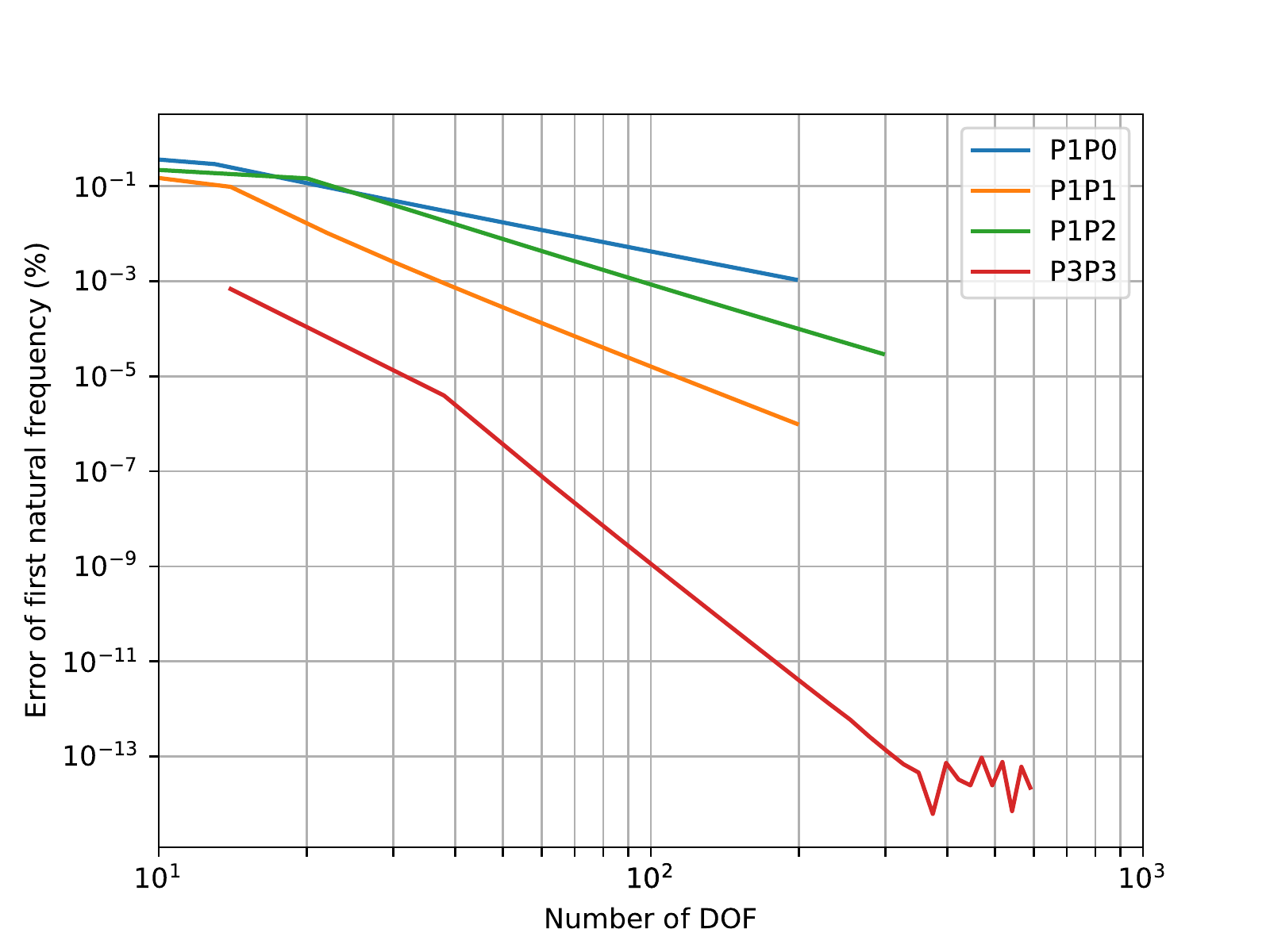}
	\end{center}
	\caption{Convergence of the first natural frequency for different polynomial interpolation of the basis functions.}
	\label{fig:convergence1d_nDOF}
\end{figure}

Recall that the variables $v_q(z,t)$, $e_q(z,t)$ and $q(z,t)$ must be discretized in $z$ with polynomials of order at least one (since they are derived once on \eqref{eq:weakform1D}). 
 Fig. \ref{fig:convergence1d_nDOF} shows the relative error of the first modal frequency for four different choices of polynomial approximations. $P_1 P_0$ stands for first order polynomial for the variables related to $q$, and order zero for the variables related to $p$. $P_1 P_1$ uses first-order polynomial for both all variables. $P_1 P_2$ uses first-order polynomial for the $q$ variables, and $p$ variables. Finally, $P_3 P_3$ uses third-order polynomial for all variables.

\FloatBarrier

Time-domain simulations for the Shallow Water Equations considering the following boundary conditions:
\begin{equation}
	\begin{split}
		e_p(0,t) = e_p(L,t) = A \cos(\omega t)\,,
	\end{split}
\end{equation}
such that these conditions represent a harmonic influx through both boundaries.

\begin{figure}[h]
	\begin{center}
	   \includegraphics[width=0.7\textwidth]{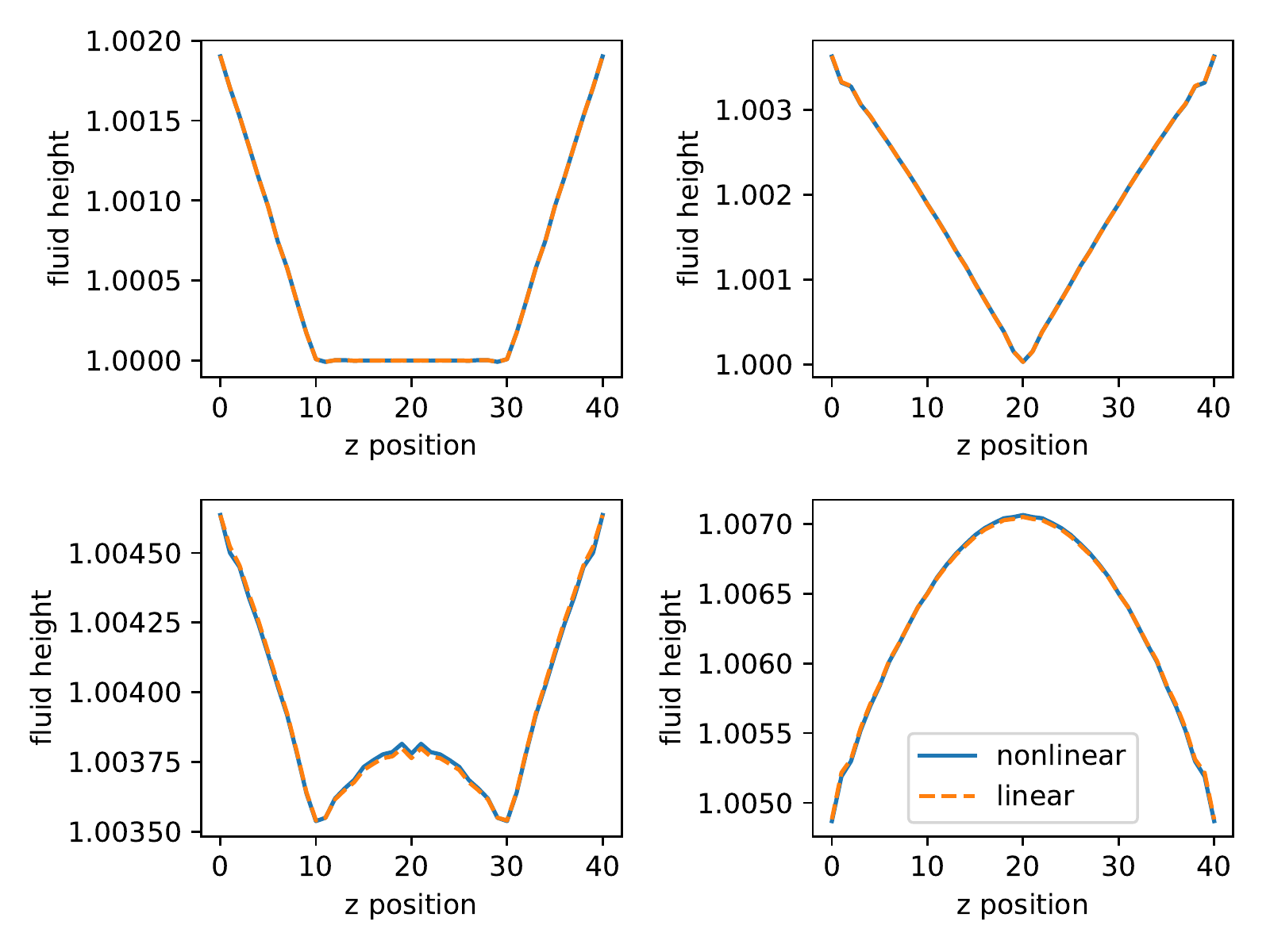}
	\end{center}
	\caption{Four snapshots of a time-domain simulation representing the fluid height as a function of horizontal position z. The previous result uses a harmonic boundary excitation with very small amplitude. The nonlinear and the linearized equations exhibts very similar results.}
	\label{fig:simu1D_Linear}
\end{figure}

\begin{figure}[h]
	\begin{center}
	   \includegraphics[width=0.7\textwidth]{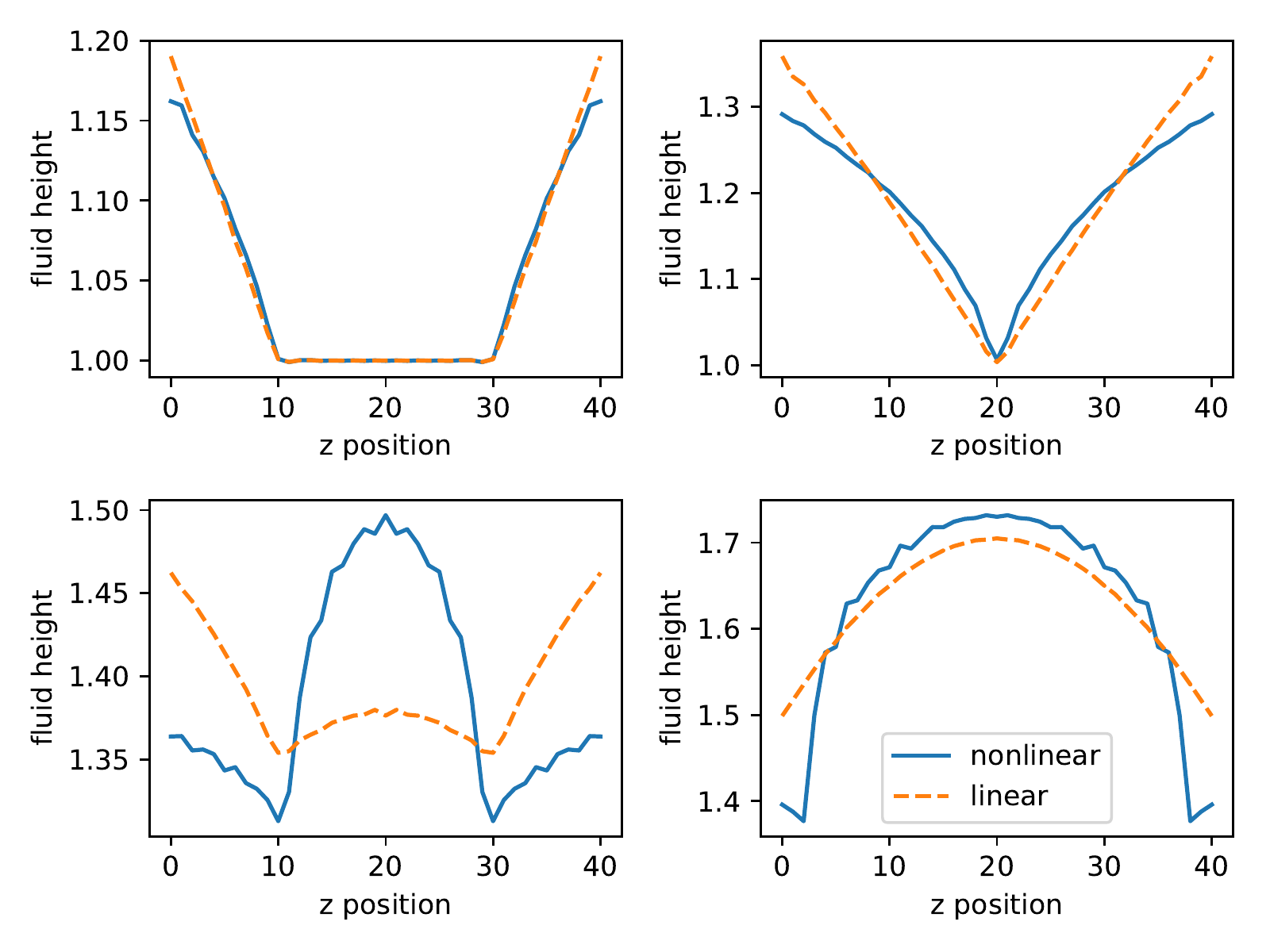}
	\end{center}
	\caption{Four snapshots of a time-domain simulation representing the fluid height as a function of horizontal position z. The previous result uses a harmonic boundary excitation with larger amplitudes, such that nonlinear phenomena now is observed.}
	\label{fig:simu1D_NonLinear}
\end{figure}

Two different amplitudes $A$ were used in the simulations. Snapshots of the simulations are presented in Figs. \ref{fig:simu1D_Linear} and \ref{fig:simu1D_NonLinear}. The first figure shows snapshots for a small amplitude value for the input inflow. In the second figure, the amplitude is multiplied by 100. Nonlinear phenomena is observed in this case.

\FloatBarrier

\FloatBarrier



\subsection{Two-dimensional Shallow Water Equations}
\label{ssNumRes2D}

The example presented in \S~\ref{sec:2DPFEM} was implemented using quadrilateral finite elements with polynomial basis functions for a square domain. 

 As we did for the 1D SWE in the previous subsection, firstly, a spectral convergence analysis of the numerical method was done.

\begin{figure}[h]
\begin{center}
   \includegraphics[width=0.7\textwidth]{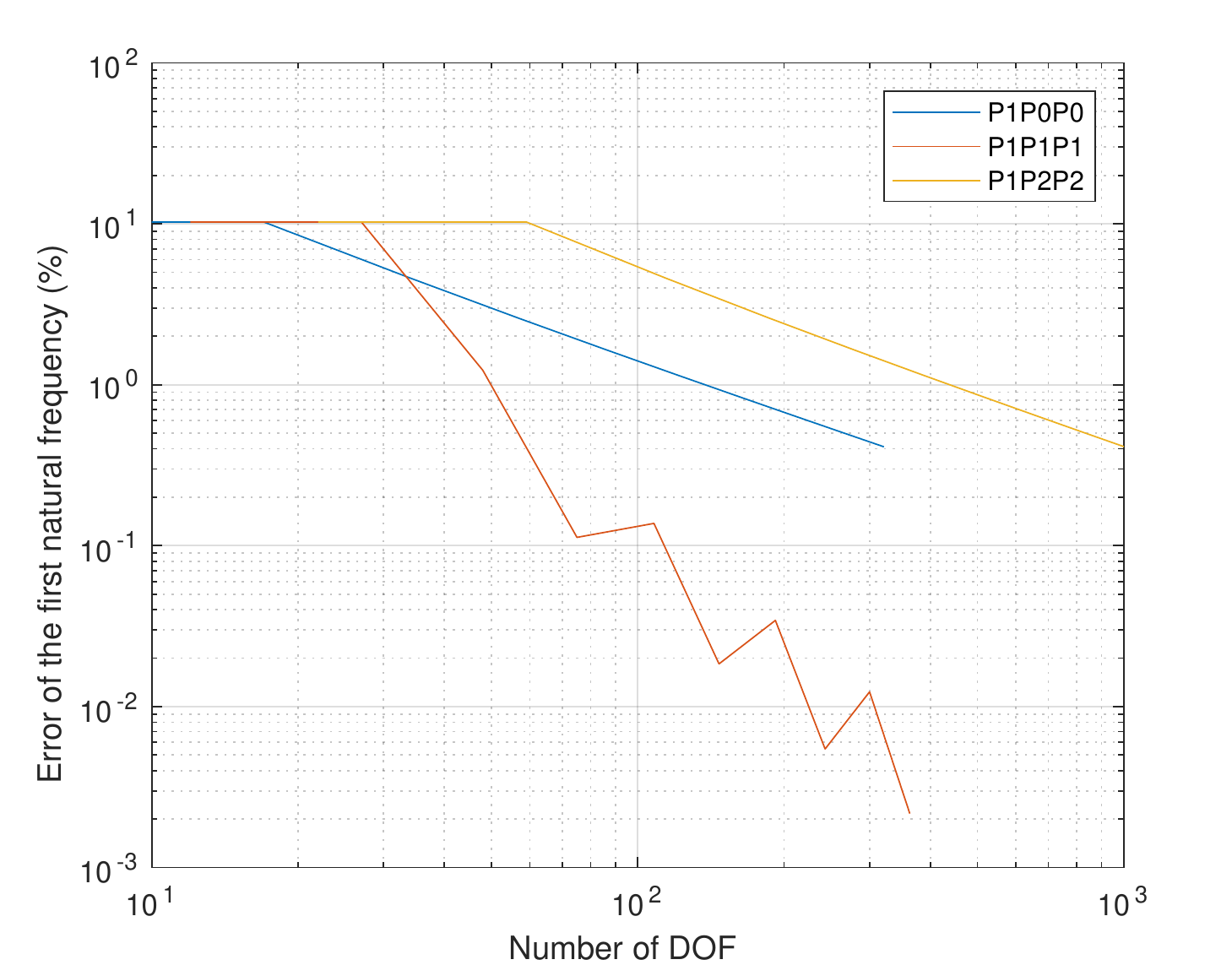}
\end{center}
\caption{Convergence of the first natural frequency of the 2D linear shallow water equations.}
\label{fig:convergence2Dmatlab}
\end{figure}

For the 2D SWE, recall that the variables with index $q$ ($\alpha_q$, $e_q$, $v_q$) must be discretized with polynomials of order at least one (since they are derived once on \eqref{eq:weakformpartitioned2D}). Fig. \ref{fig:convergence2Dmatlab} shows the relative error of the first modal frequency for three different choices of polynomial approximations. $P_i P_j P_j$ stands for "i-th" order polynomial for the variables with index $q$, and "j-th" order for the two components of the vector variables of index $p$ ($\pmb{\alpha_p}$, $\pmb{e_p}$ and $\pmb{v}_p$). 

\FloatBarrier

\begin{figure}[h]
	\begin{minipage}{\linewidth}\centering
		\subcaptionbox{t = 0.4 s}{\includegraphics[width=0.49\textwidth]{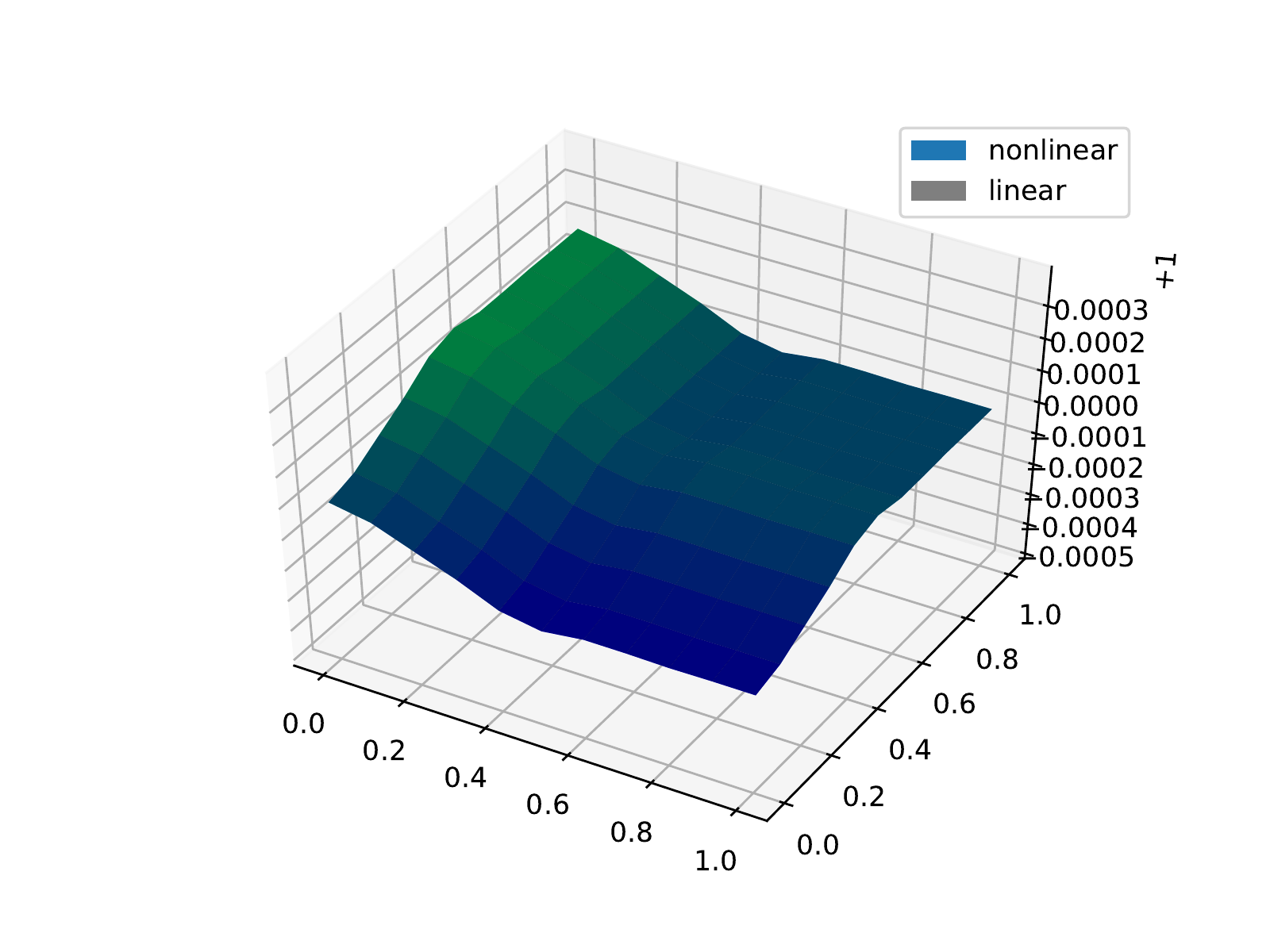}
		}
		\subcaptionbox{t = 0.8 s}{\includegraphics[width=0.49\textwidth]{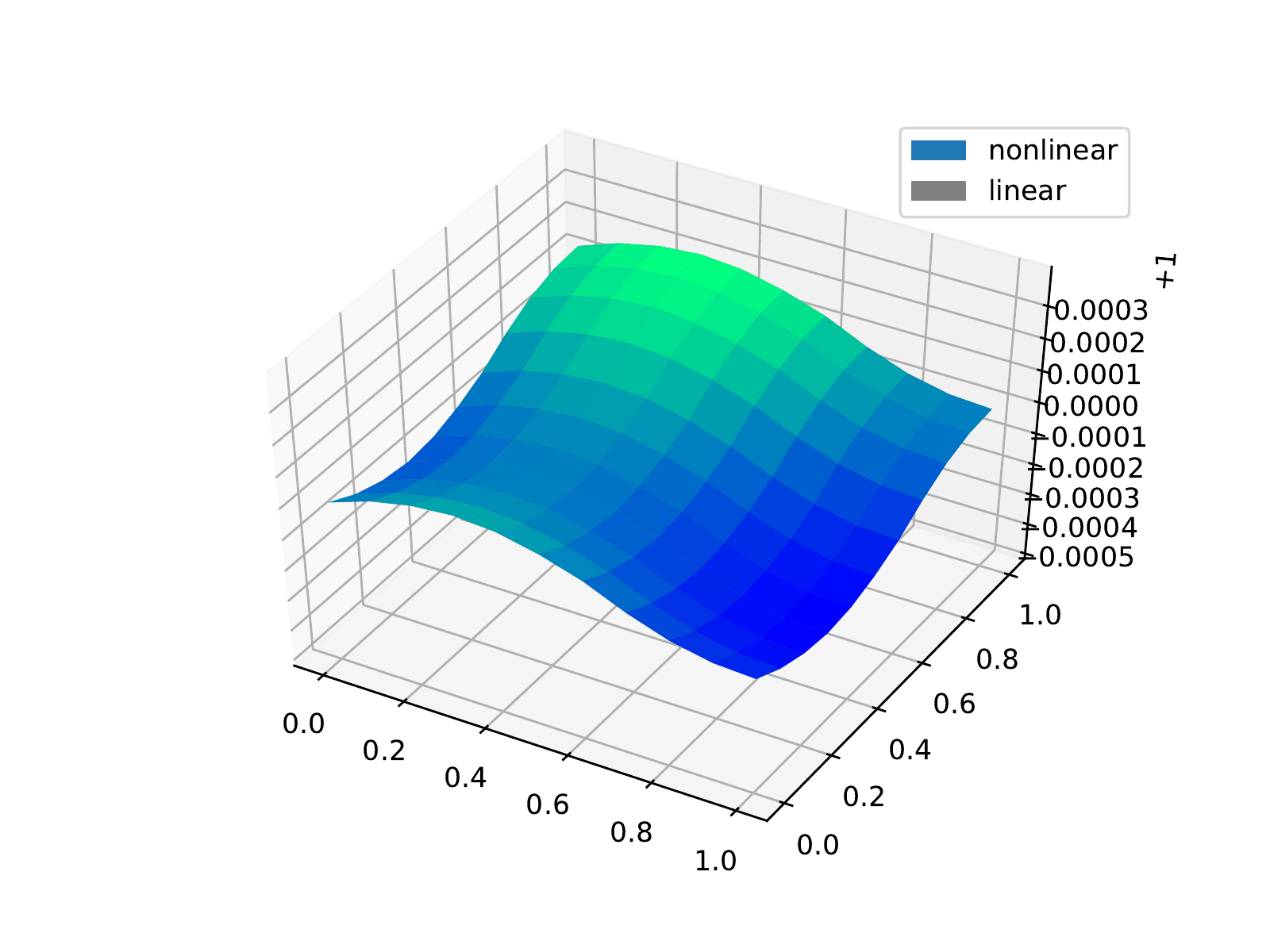}
		}
		\subcaptionbox{t = 1.2 s}{\includegraphics[width=0.49\textwidth]{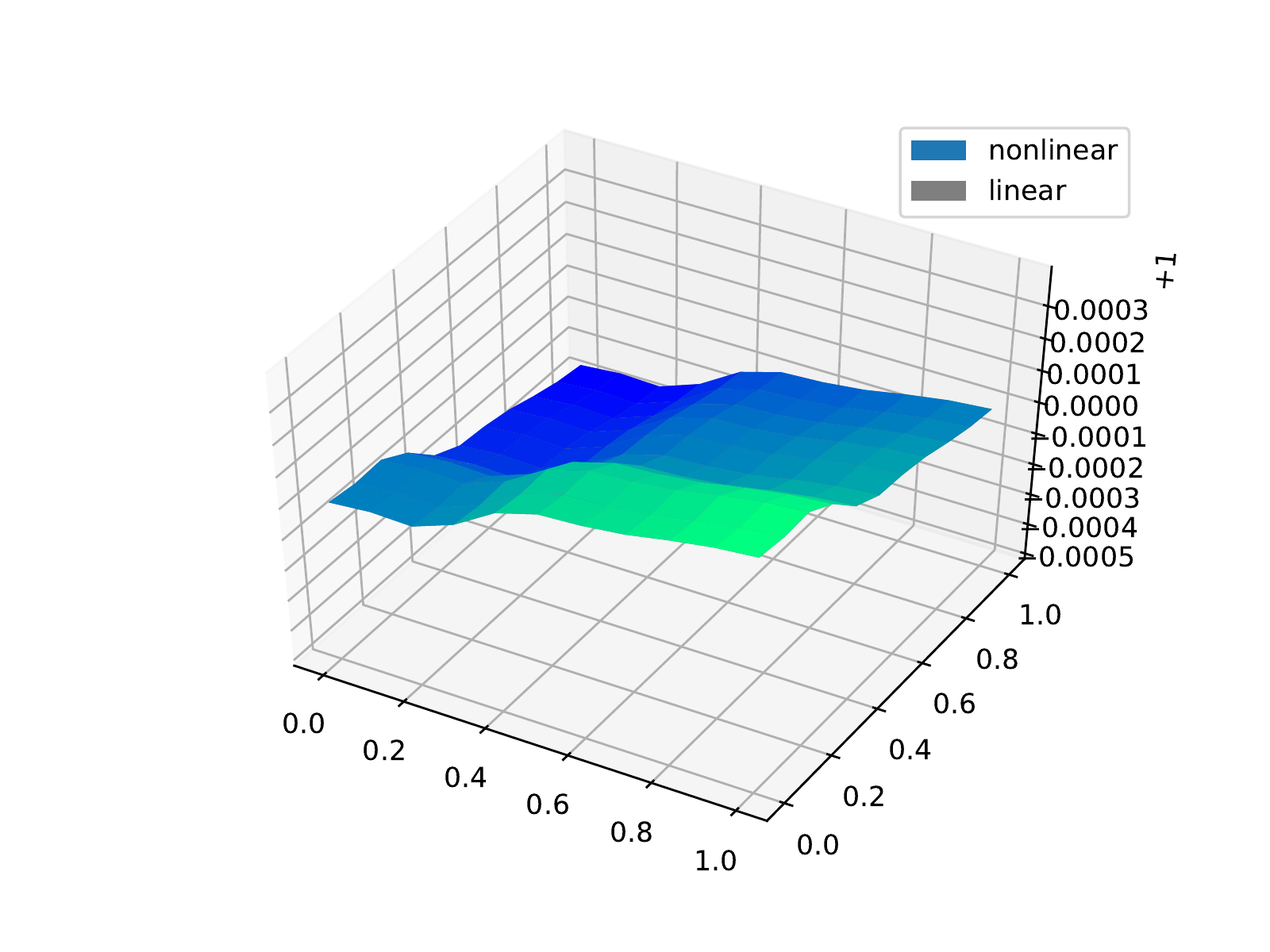}
		}
		\subcaptionbox{t = 1.6 s}{\includegraphics[width=0.49\textwidth]{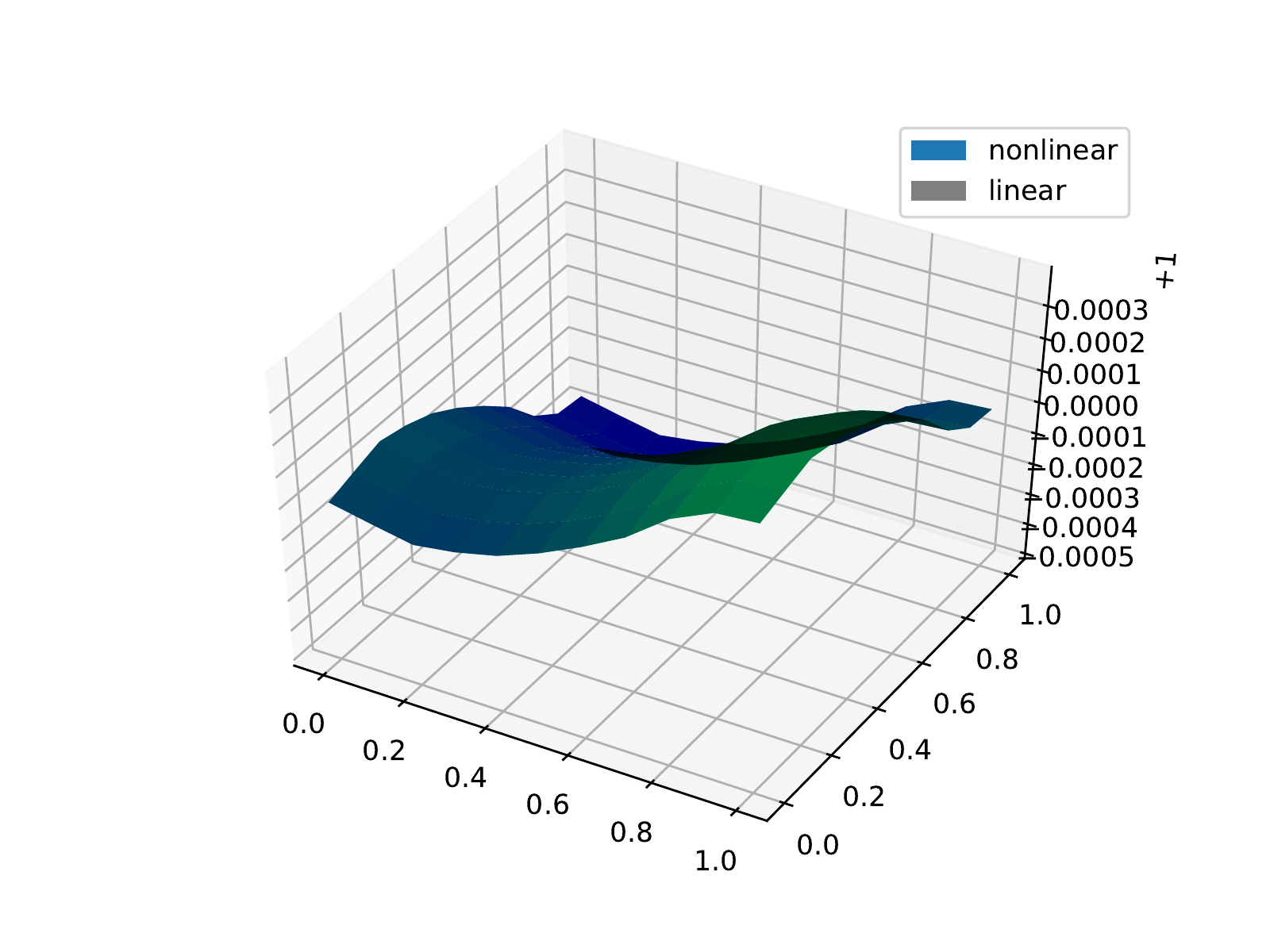}
		}
		\caption{Snapshots of simulation for a harmonic inflow excitation at two of the boundaries of the domain. The variable $\alpha_q$ (fluid height) is shown. Here, small inputs are considered, such that the nonlinear and linearized time-responses are almost equivalent.}
	\label{fig:simulation2D_borderinflow_exc_linear}
	\end{minipage}
	\end{figure}

\begin{figure}[h]
		\centering
		\includegraphics[width=0.8\textwidth]{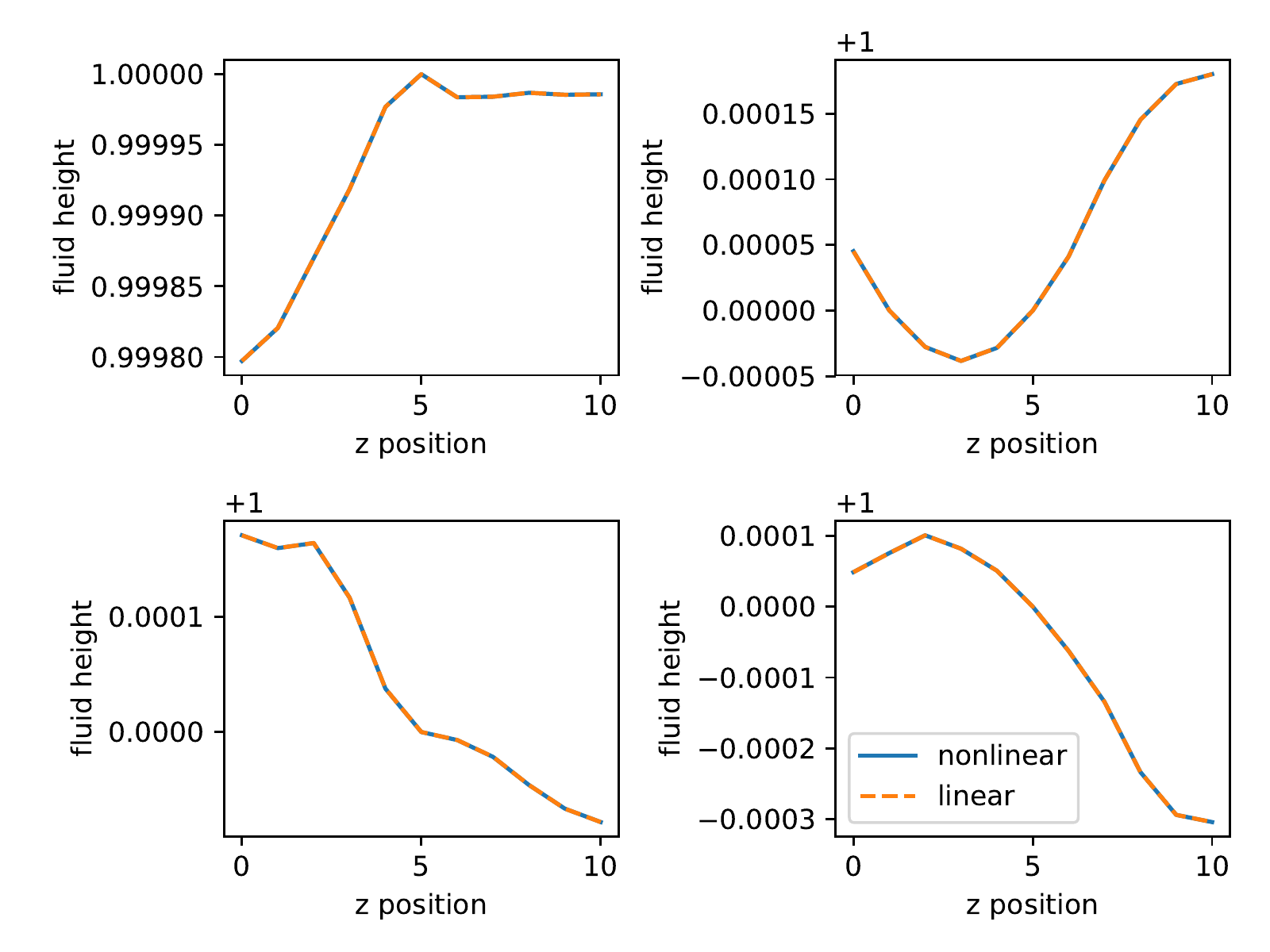}
		\caption{Snapshots of simulation for a harmonic inflow excitation at two of the boundaries of the domain. The variable $\alpha_q$ (fluid height) is shown along a cross-section in the middle of the domain. Here, small inputs are considered, such that the nonlinear and linearized time-responses are almost equivalent.}
		\label{fig:simulation2D_borderinflow_exc_linear_border}
	\end{figure}

\begin{figure}[h]
\begin{minipage}{\linewidth}\centering
	\subcaptionbox{t = 0.4 s}{\includegraphics[width=0.49\textwidth]{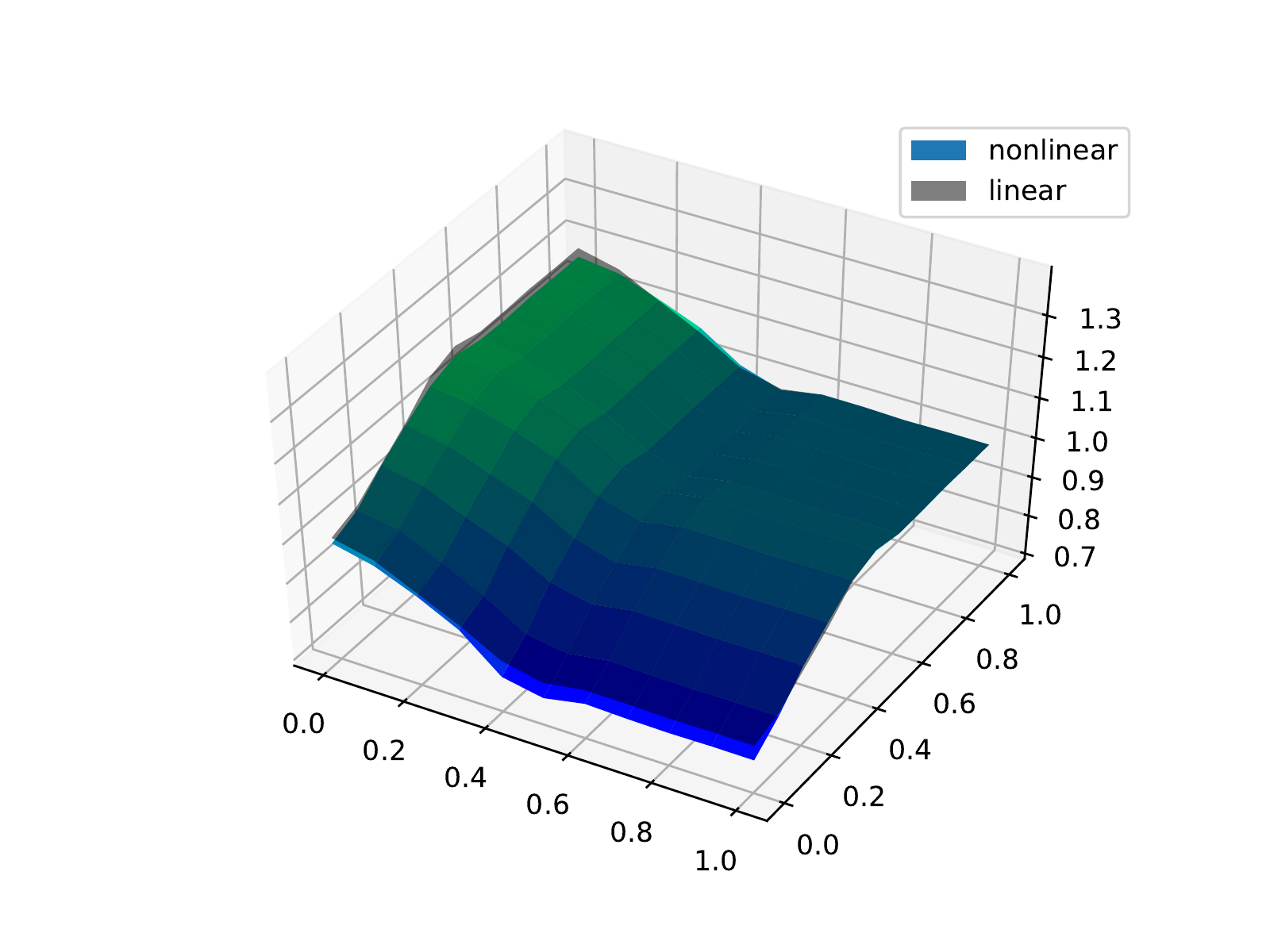}
	}
	\subcaptionbox{t = 0.8 s}{\includegraphics[width=0.49\textwidth]{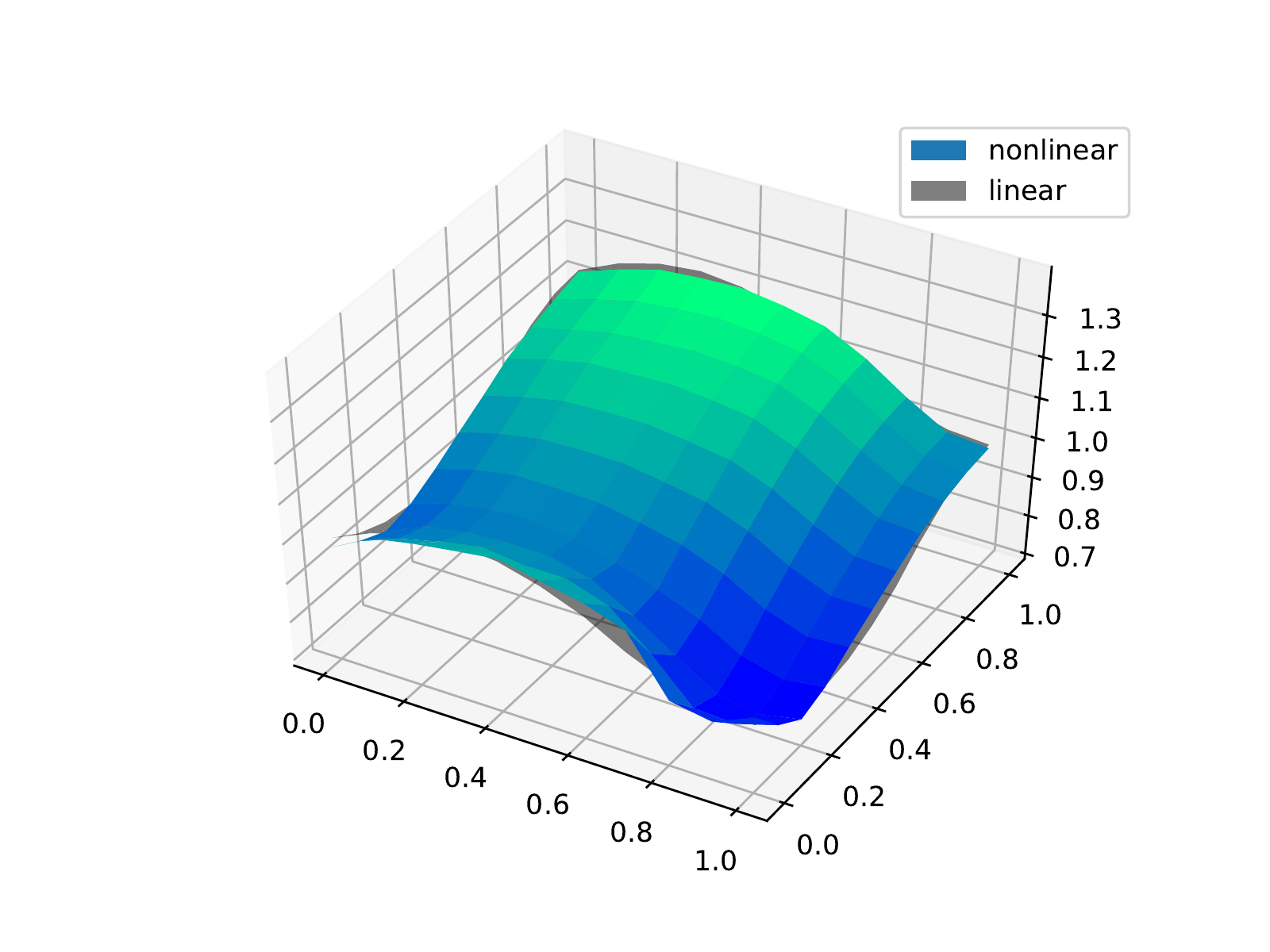}
	}
	\subcaptionbox{t = 1.2 s}{\includegraphics[width=0.49\textwidth]{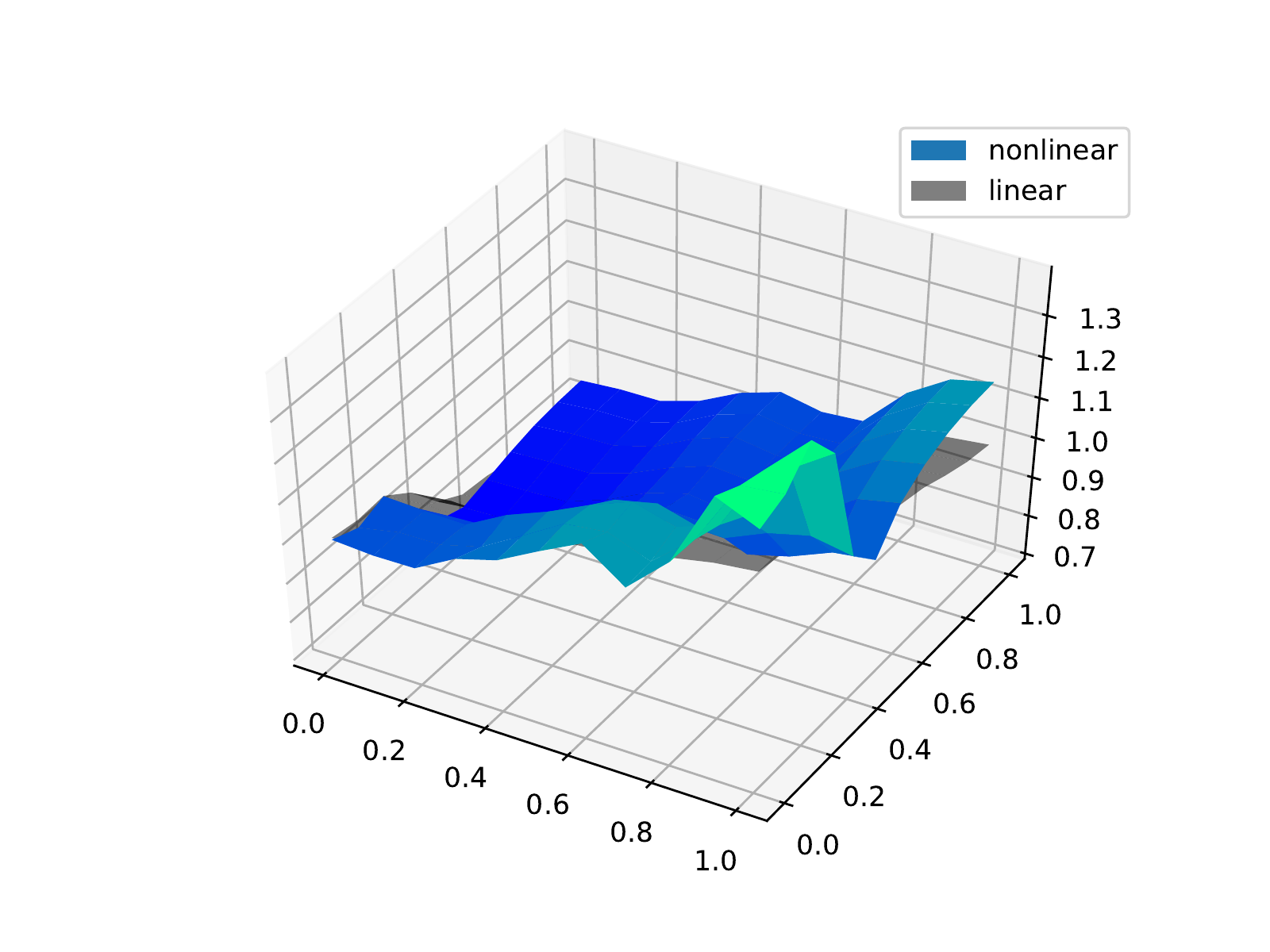}
	}
	\subcaptionbox{t = 1.6 s}{\includegraphics[width=0.49\textwidth]{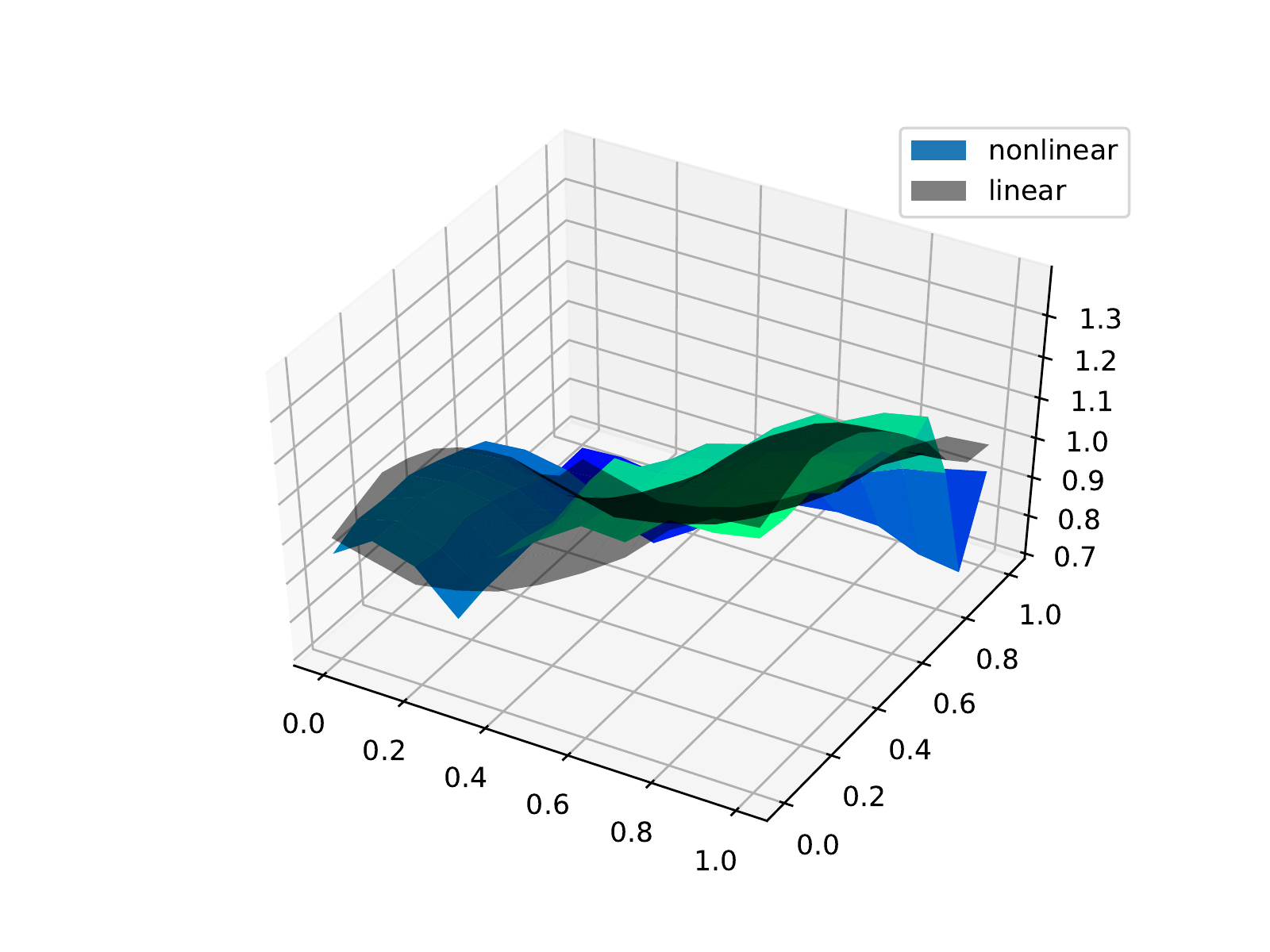}
	}
	\caption{Snapshots of simulation for a harmoninc inflow excitation at two of the boundaries of the domain. The variable $\alpha_q$ (fluid height) is shown. Differences between the nonlinear and linearized time-responses are now observed.}
	\label{fig:simulation2D_borderinflow_exc_nonlinear}
\end{minipage}
\end{figure}

\begin{figure}[h]
	\centering
	\includegraphics[width=0.8\textwidth]{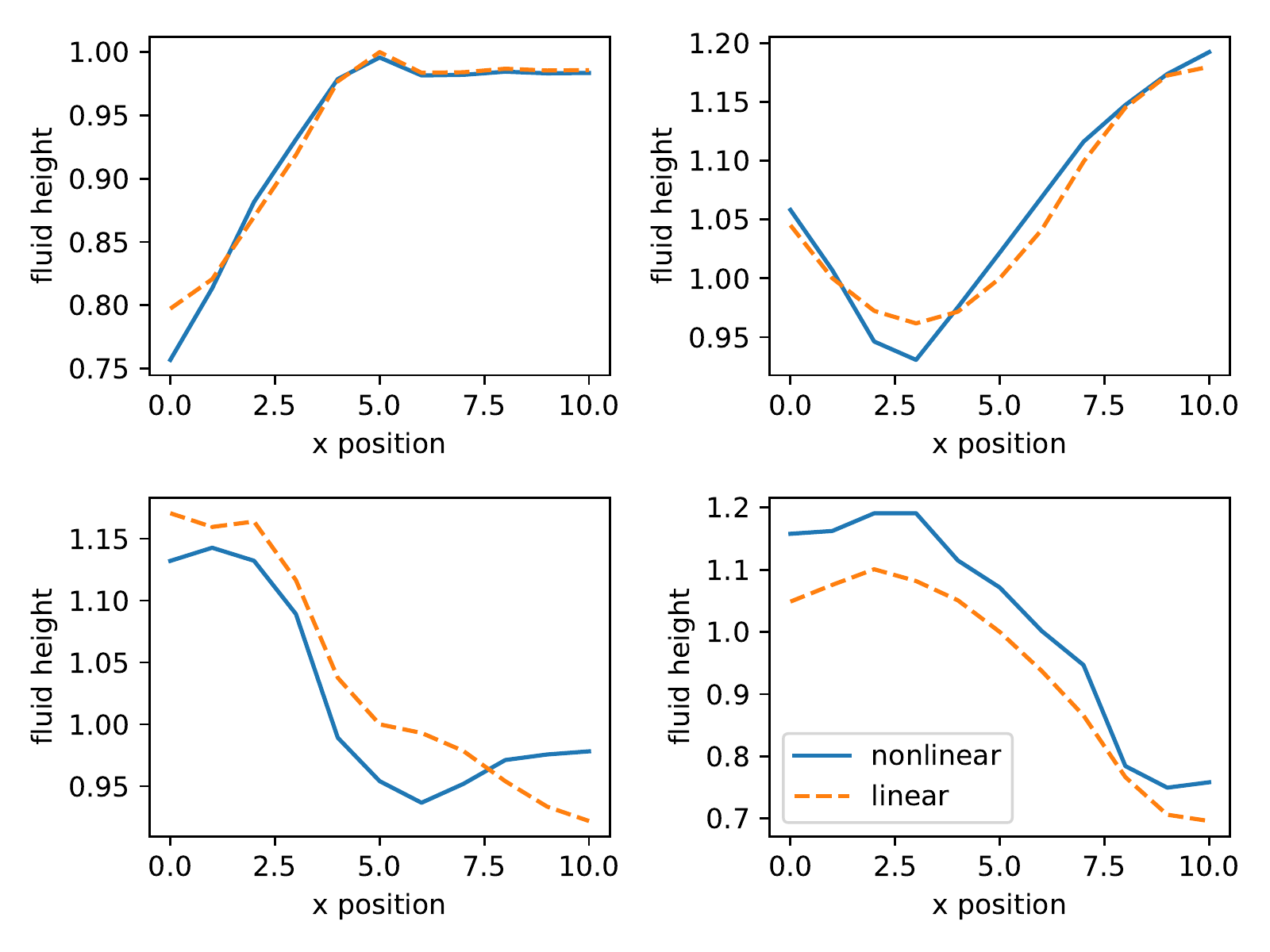}
	\caption{Snapshots of simulation for a harmoninc inflow excitation at two of the boundaries of the domain. The variable $\alpha_q$ (fluid height) is shown along a cross-section in the middle of the domain. Differences between the nonlinear and linearized time-responses are observed.}
	\label{fig:simulation2D_borderinflow_exc_nonlinear_border}
\end{figure}

Numerical simulations were performed using the discretized system under boundary-port excitation.  The following boundary conditions were considered:
\begin{equation}
	u_\partial(x,y,t) = \left\{
	\begin{split}
		A \sin(\pi t) \,\,, &[x,y] \in \partial\Omega_{up}\,,\\
		-A \sin(\pi t) \,\,, &[x,y] \in \partial\Omega_{left}\,,\\
		0 \,\,, &[x,y] \in \partial\Omega_{down} \cup \partial\Omega_{right}\,.
	\end{split}
	t \leq 1 s,
	\right.
\end{equation}
and $u_\partial(x,y,t) = 0, [x,y] \in \partial\Omega, t > 1 s$.
The boundary is split in four sides: $\partial\Omega = \partial\Omega_{up} \cup \partial\Omega_{left} \cup \partial\Omega_{down} \cup \partial\Omega_{right}$.
These conditions impose a harmonic inflow on one side of the boundary and the opposite condition on the other side. Simulations for two different values of amplitude $A$ are presented. First, snapshots for small amplitudes are presented in Figs. \ref{fig:simulation2D_borderinflow_exc_linear} and \ref{fig:simulation2D_borderinflow_exc_linear_border}.
Secondly, snapshots for large amplitudes are presented in Figs. \ref{fig:simulation2D_borderinflow_exc_nonlinear} and \ref{fig:simulation2D_borderinflow_exc_nonlinear_border}.

\begin{figure}[h]
\begin{minipage}{\linewidth}\centering
		\subcaptionbox{Hamiltonian}{\includegraphics[width=0.49\textwidth]{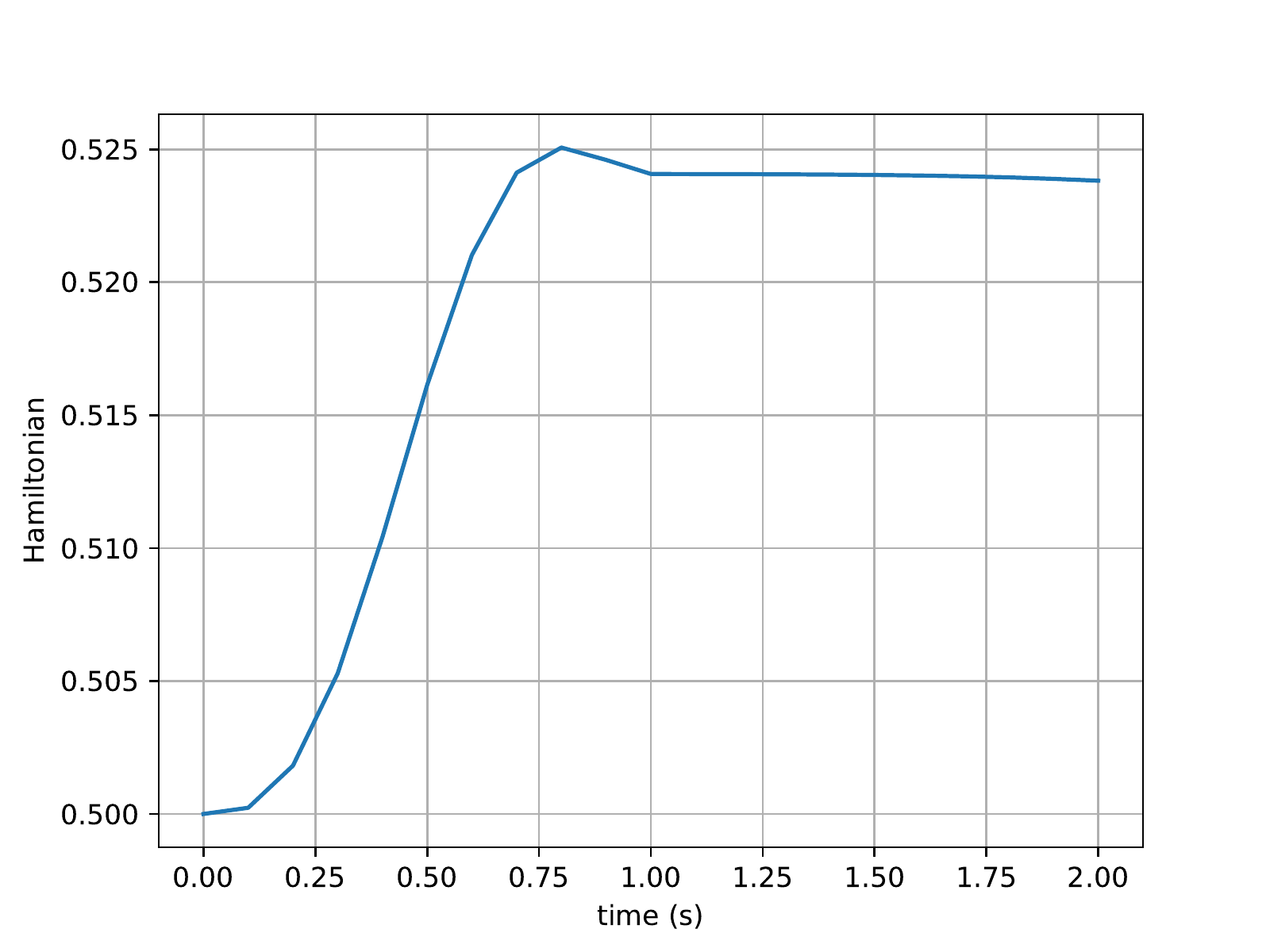}
		}
		\subcaptionbox{Total volume}{\includegraphics[width=0.49\textwidth]{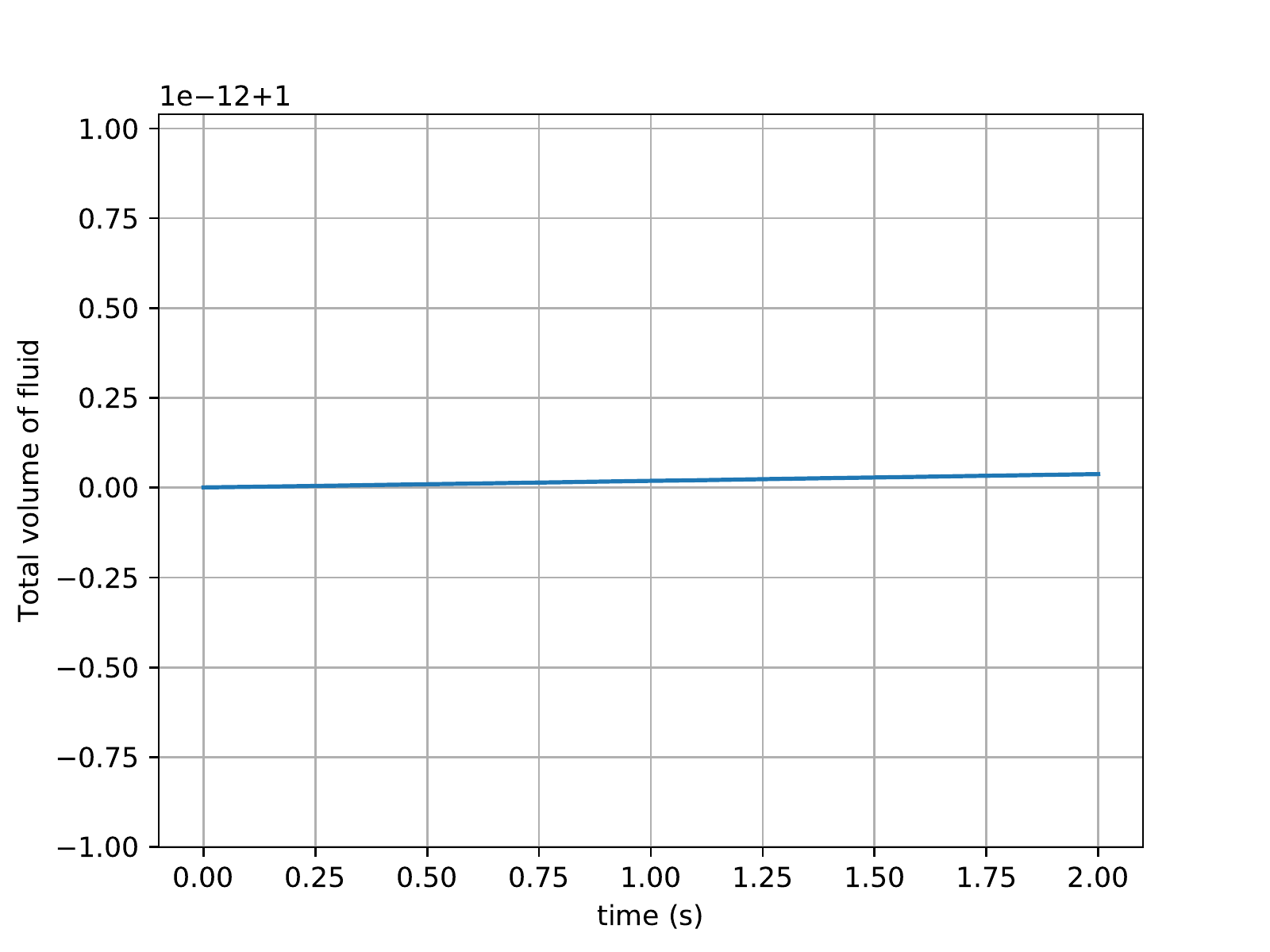}
		}
\caption{a) Hamiltonian as a function of time for the simulation for a harmoninc inflow excitation at two of the boundaries of the domain. The Hamiltonian only changes during the first second, while an external excitation is applied. b) Total volume of fluid as a function of time for the simulation for a harmonic inflow excitation at two of the boundaries of the domain. The total volume is constant along all the simulation.}
\label{fig:HamiltonianVolume2Dsimu}
\end{minipage}
\end{figure}

Figure \ref{fig:HamiltonianVolume2Dsimu} shows how the approximated Hamiltonian \eqref{eq:2DSWEconstrel} and total volume of fluid ($V(t) = \int_\Omega \alpha_1^{ap}(x,y,t) \dOmega$) change with time. As expected, the Hamiltonian only changes during the first second of simulation, while the system is excited through the boundary ports. After that, since $\dot{H}_d = \pmb{u}_\partial^T \pmb{y}_\partial = 0$, the Hamiltonian is constant. Finally, the total volume is kept nearly constant and the changes are only due to numerical precision (of order $10^{-14}$).

\FloatBarrier

\section{Extensions}
\label{sec:extensions}
In this final section, the extension of  PFEM to more general configurations is addressed, in order to illustrate the flexibility of the proposed numerical method: in \S~\ref{ss-polar} the specific expression in polar coordinates is investigated, in \S~\ref{ss-inhomogeneous} the case of heterogeneous medium with variable coefficients is presented, and in \S~\ref{ss-highorder} the method is applied to second-order differential operators, such as those involved in the port-Hamiltonian formulation of the Euler-Bernoulli beam in 1D.

\subsection{Polar coordinates}
\label{ss-polar}
The goal of this subsection is to prove the applicability of PFEM in 2D,  when the chosen coordinate system is not Cartesian: polar coordinates are presented.
In order to avoid unnecessary technicalities, the geometry of a disc has been chosen to illustrate this extension.

\paragraph{The 2D Shallow Water Equations in a disc, as a Port-Hamiltonian system in polar coordinates}
Let us consider the disc $\Omega=D_R$ of radius $R>0$ with boundary $\partial \Omega=C_R$, the circle of radius $R$. Polar coordinates $r$ and $\theta$ will be used.
In vector calculus, the  2-form $\alpha^q= h\,(r \dr \wedge \dth)$ and the 1-form $\alpha^p= \rho \, (u^r\, \dr + u^{\theta}\, r\dth)$ are represented by the scalar function $\alpha_q= h$, and the vector function $\pmb{{\alpha}_p}:=\rho\,[u^r(t,r,\theta),u^{\theta}(t,r,\theta)]^T$ respectively.
The Hamiltonian reads
\begin{eqnarray}
  H&=& \frac{1}{2}\int_{D_R} [\rho g\, h^2 + \rho h\,((u^r)^2+(u^{\theta})^2)]  \, r\,\dr\dth\,, \\
  &=& \int_{D_R} [\frac{1}{2}\rho g\, \alpha_q^2 + \frac{1}{2\rho}  \alpha_q |\pmb{{\alpha}_p}|^2]\,  r\,\dr\dth\,.
  \end{eqnarray}
The effort or co-energy variables can be computed as
$e_q :=\delta_q H=\rho g\, \alpha_q +\frac{1}{2\rho}|\pmb{{\alpha}_p}|^2$, and 
$\pmb{e_p} :=\delta_{\pmb{p}} H=\frac{1}{\rho}  \alpha_q \pmb{{\alpha}_p} =h\,[u^r(t,r,\theta), u^{\theta}(t,r,\theta)]^T$ (which is associated to the 1-form  $e^p =h \left( u^{\theta} \dr - u^r r \dth \right)$.
Finally, the 2-form $\mbox{d}\,e^p$ translates into the scalar $\mbox{div}(\pmb{e_p}):= r^{-1}\,\partial_{r}(r\,e_p^r)   + r^{-1}\,\partial_{\theta} e_p^{\theta}$, and the 1-form $\mbox{d}\,e^q$ translates into the vector $\pmb{grad}(e_q):= [\partial_r e_q, \, r^{-1}\,\partial_{\theta}  e_q]^T$.
With these notations and definitions, we get the same system as (\ref{eq:2DSWECartCoord}) for the strong form of the pHs, namely:
\begin{equation}
	\label{eq:2DSWECartCoordPolar}
	\left[ \begin{array}{c} \dot{h} \\ \rho  \left[ \begin{array}{c} \dot{u^r} \\ \dot{u^{\theta}} \end{array} \right] \end{array} \right]  = \left[ \begin{array}{cc} 0 & - \mbox{div} \\ -\pmb{grad} & 0 \end{array} \right] \left[ \begin{array}{c}  \rho(g h + \frac{(u^r)^2+(u^{\theta})^2}{2}) \\ h \left[ \begin{array}{c} {u^r} \\ {u^{\theta}} \end{array} \right] \end{array} \right]\,, 
\end{equation}
with boundary control $u_{\partial}(\theta,t):=-\pmb{e_p} \cdot \pmb{n}=-e_p^r(R, \theta, t)$ and collocated boundary observation $y_{\partial}(\theta,t):=e_q(R, \theta, t)$ at the boundary $\partial \Omega=C_R$.
Let us conclude with the energy balance for this system:
\begin{equation}
	\label{eq:Hamiltonianpolar}
        \frac{d}{dt} \frac{1}{2}\int_{D_R} [\rho g\, h^2 + \rho h\,((u^r)^2+(u^{\theta})^2)] \,r\, \dr \dth=\int_{C_R} u_{\partial}(\theta,t)\,y_{\partial}(\theta,t)\,R\,\dth\,.
\end{equation}

\paragraph{The Partitioned Finite Element Method directly applies to the  port-Hamiltonian system in polar coordinates}
Now, let us first  rewrite the weak form (\ref{eq:weakformpartitioned2D}) with test functions $v_q$ and $\pmb{v_p}$.
\begin{equation}
	\label{eq:weakformpartitioned2Dpolar}
	\begin{split}
		\int_{D_R} {v_q {\dot{\alpha}_q}}\,r\, \dr \dth & = \int_{D_R} { \left( \pmb{\nabla}v_q \right) \cdot \pmb{e}_{\pmb{p}}}\,r\, \dr \dth - \int_{C_R} v_q \, \pmb{n} \cdot \pmb{e}_{\pmb{p}} \,R\,\dth\,,\\
		\int_{D_R}  {\pmb{v_p} \cdot \pmb{\dot{\alpha}_p}}\,r\, \dr \dth & = - \int_{D_R} {\pmb{v_p} \cdot  \pmb{\nabla} e_q}\,r\, \dr \dth \,.
	\end{split}
\end{equation}

Let us approximate the scalar energy variables  $\alpha_q(r,\theta,t)$ using the following basis with $N_q$ elements:
\begin{equation}
	\alpha_q(r,\theta,t) \approx \alpha_q^{ap}(r,\theta,t) := \sum_{i=1}^{N_q} \phi^i_q(r,\theta) \alpha_q^i(t) = \pmb{\phi}_q(r,\theta)^T \pmb{\alpha}_q(t)\,.
\end{equation}
The variables $e_q$ and $v_q$ are also approximated using $\pmb{\phi}_q(r,\theta)$.

Similarly, the vectorial energy variable $\pmb{\alpha_p}$ is approximated as:
\begin{equation}
	\pmb{\alpha_p}(r,\theta,t) \approx \pmb{\alpha_p}^{ap}(r,\theta,t) := \sum_{k=1}^{N_p} \pmb{\phi_p}^k(r,\theta) \alpha_p^k(t)  = \pmb{\Phi}_p(r,\theta)^T \pmb{\alpha}_p(t)\,,
\end{equation}
where ${\pmb{\phi_p}^k(r,\theta) = \left[\begin{matrix} \phi_p^{r,k} (r,\theta) \\ \phi_p^{\theta,k} (r,\theta) \end{matrix}\right]}$ represents a 2D-vectorial basis function and, consequently, $\pmb{\Phi_p}(r,\theta)$ is an $N_p \times 2$ matrix. Remark~\ref{rem6} does apply here also.
Furthermore, $\pmb{e_p}$ and $\pmb{v_p}$ are also approximated using $\pmb{\Phi}_p(r,\theta)$. 

Finally, the boundary input, localized on the circle of radius $r=R$ can be discretized using any one-dimensional set of basis functions, say $\pmb{\psi}=[\psi^m]$, provided $2\pi$-periodicity is ensured (trigonometric polynomials are a fair trial  approximation basis, see e.g. \cite[Chapter~18]{Boyd2001}):
 \begin{equation}
	u_\partial(\theta,t) \approx u^{ap}_\partial(\theta,t) := \sum_{m=1}^{N_\partial} \psi^m(\theta) u_{\partial}^m(t) = \pmb{\psi}(\theta)^T \pmb{u}_{\partial}(t)\,.
 \end{equation}

 Introducing the notation $\partial_r \pmb{\phi}_q:=[\partial_r {\phi}_q^i ]$ and $\partial_{\theta} \pmb{\phi}_q:=[\partial_{\theta} {\phi}_q^i ]$ for the matrices of partial derivatives of the functions ${\phi}_q^i$, we define matrix
 \begin{equation}
 D:=\int_{D_R}{ \left[ \begin{matrix} \partial_r\pmb{\phi}_{q} & r^{-1}\,\partial_{\theta}\pmb{\phi}_{q} \end{matrix} \right] \pmb{\Phi}_p^T r\,\dr\dth}
=\int_{D_R}{ \left[ \begin{matrix} r\,\partial_r\pmb{\phi}_{q} & \partial_{\theta}\pmb{\phi}_{q} \end{matrix} \right] \pmb{\Phi}_p^T \dr\dth}\,,
\end{equation}
where the apparent singularity at $r=0$ has been removed.
Then, with classical mass matrices $M_q:=\int_{D_R} {\pmb{\phi}_q \pmb{\phi}_q^T r\,\dr\dth}$, $M_p:=\int_{D_R}{\pmb{\Phi}_p \pmb{\Phi}_p^T r\,\dr\dth}$,  
together with the control matrix $B:= \int_{C_R}{\pmb{\phi}_q(R,\theta)}\,\pmb{\psi}^T(\theta) \, R\,\dth$, the finite-dimensional equations become:
\begin{equation}
	\label{eq:finitedimensionaleqs}
	\begin{split}
		{M_q}\, \dot{\pmb{\alpha}}_q  = & {D}\, \pmb{e}_p + {B} \,\pmb{u}_\partial(t) \,,\\
		{M_p}\, \dot{\pmb{\alpha}}_p = &  - {D^T}\, \pmb{e}_q \,,
	\end{split}
 \end{equation}
where $M_q$ and $M_p$ are square matrices (of size $N_q \times N_q$ and $N_p \times N_p$, respectively). $D$ is an $N_q \times N_p$ matrix and $B$ is an $N_q \times N_\partial$ matrix.

Defining $\pmb{y}_\partial(t)$, the output conjugated to the input $\pmb{u}_\partial(t)$ as:
 \begin{equation}
	\pmb{y}_\partial(t) :={M_{\psi}}^{-1} B^T \pmb{e}_q(t)\,,
 \end{equation}
with boundary mass matrix $M_{\psi}:=\int_{C_R}{\pmb{\psi} \pmb{\psi}^T R\,\dth}$,  the approximated system can be written using the finite-dimensional Dirac structure representation given by (\ref{eq:finitedim_2dwaveeq}), and as found in remark~\ref{rem8}, the global energy balance reads
\begin{equation}
	\label{eq:finitedimHamiltonianpolar}
\frac{d}{dt} \frac{1}{2}(\pmb{\alpha}_q^T {M_q} \pmb{\alpha}_q + \pmb{\alpha}_p^T {M_p} \pmb{\alpha}_p) = \pmb{y}_\partial^T {M_{\psi}} \pmb{u}_\partial\,,
\end{equation}
which mimicks that at the continuous level, namely (\ref{eq:Hamiltonianpolar}).
 
\subsection{Heterogeneous case with variable coefficients}
\label{ss-inhomogeneous}
The goal of this subsection is to prove the applicability of PFEM when the coefficients are space varying. In order to avoid unnecessary technicalities, the choice has been made to tackle the 1D model, first derived in \S~\ref{sec:introexample} as introductory example. Another fully worked out example can be found in~\cite{SerhaniCPDE2019} on the anisotropic heterogeneous wave equation in 2D.

\paragraph{The variable-coefficient physical model as a port-Hamiltonian system} Let us consider the Shallow Water equation in a water channel with a space-varying cross section, i.e. with $z \mapsto b(z)$ the width of the channel, it is easy to understand that the energy given by~(\ref{eq:SWE1H}) remains unchanged, with function $b(z)$ instead of coefficient $b$, but then then the balance equations ~(\ref{eq:SWE}) must be modified as follows:
\begin{equation}
	\label{eq:SWEvariable}
	\begin{split}
		\frac{\partial}{\partial t}(b\,h) = - \frac{\partial}{\partial z}\left(b h u\right) \,,\\
		\frac{\partial}{\partial t}u = - \frac{\partial}{\partial z}\left(\frac{u^2}{2} + gh\right) \,.
	\end{split}
\end{equation}
With the appropriate choice of energy variables $q:= b\,h$, and $p:=\rho\,b\,u$, the Hamiltonian (\ref{eq:SWE1H}) now reads:
\begin{equation}
	\label{eq:waveequation_Hamvariable}
	H\left( q(z,t), p(z,t)\right) = \frac{1}{2}\int_{[0,L]} \left( \frac{q p^2}{\rho\,b^2} + \, \frac{\rho g}{b} q^2 \right) \dz \,.
 \end{equation}
The co-energy variables are found to be $e_q:=\delta_q {\cal H}=\rho \left( \frac{u^2}{2} + g h \right)$,  and $e_p:=\delta_p {\cal H} = h\,u$; thus system~(\ref{eq:SWEvariable}) becomes in compact form:
\begin{equation}
	\label{eq:1Dstrongformvariable}
	\begin{split}
		\dot{q}(z,t) & = - \frac{\partial}{\partial z}\left[ b(z)\,e_p (z,t) \right] \,,\\
		\dot{p}(z,t) & = - b(z)\,\frac{\partial}{\partial z} e_q (z,t)\,,
	\end{split}
 \end{equation}
to be compared with (\ref{eq:1Dstrongform}) in the uniform case.
Hence the new interconnection operator ${\cal J}_b$ reads:
$$
{\cal J}_b:=\begin{bmatrix} 0 & -{\partial_z}\left[b(z)\,. \right] \\
- b(z)\,{\partial_z} & 0
\end{bmatrix}\,.
$$
Since $\int_{[0,L]} \varphi\,\partial_z (b(z)\,\psi)\, \dz = - \int_{[0,L]} b(z)(\partial_z\varphi)\,\psi\, \dz$ for smooth scalar functions $\varphi$ and $\psi$, a straightforward computation shows that $({\cal J}_b\pmb{u}\,,\pmb{v}) = - (\pmb{u}\,,{\cal J}_b \pmb{v})$ for smooth vector-valued functions $\pmb{u}$ and $\pmb{v}$ vanishing at the ends of the interval, and with the standard scalar product in $L^2 \times L^2$. Hence, the unbounded matrix-valued differential operator ${\cal J}_b$ proves skew-symmetric in $L^2 \times L^2$.

\begin{remark}
The port-Hamiltonian formulation \eqref{eq:1Dstrongformvariable} makes use of a non canonical Stokes-Dirac structure with non uniform coefficents. Instead, one could use the intrinsic geometric formulation presented in subsection \ref{EDFSDS}, together with constitutive equations in covariant form for the SWE in subsection \ref{subsec:wave2d} and make use of the same energy state variables as in the uniform case, that is:
\begin{eqnarray}
\alpha^q&:=&h(z,t) b(z) \, dz \nonumber \\
\alpha^p&:=&\rho u(z,t) \, dz \nonumber
\end{eqnarray}
In this approach, information about the geometry - here the space-varying reach width - will be embedded in the constitutive equation, that is in the Hamiltonian functional, through the definition of the appropriate Hodge star operator. More precisely, using the non uniform metric derived from the duality product $({u},{v}) = \int_{[0,L}{u(z)v(z)\, b(z)dz}$, one gets
\begin{eqnarray}
\star\left( \alpha^q \right) &:=&h(z,t) \nonumber \\
\star\left( \alpha^p \right) &:=&\rho u(z,t) b(z) \nonumber
\end{eqnarray}
Therefore the Hamiltonian functional \eqref{eq:HSW2Ddef} reads
\begin{eqnarray}
H(\alpha^q,\alpha^p)&:= &\int_\Omega{ \frac{\rho g\,(\star{\alpha^q}) \alpha^q}{2} + \frac{\star{\alpha^q}\left( \alpha^p \wedge \star{\alpha^p} \right)}{2 \rho} } \nonumber \\
&=& \frac{\rho}{2} \int_{[0,L]}{ gb(z)h^2(z,t) + b(z)h(z,t)u^2(z,t) } \, dz \nonumber
\end{eqnarray}
and the corresponding co-energy variables defined in \eqref{eq:eqepSW2Ddef} read
\begin{eqnarray}
		e_q  & = & \delta_q H = \rho g (\star{\alpha^q}) + 
\frac{1}{2\rho} \star \left( (\star{\alpha^p}) \wedge \alpha^p \right) = \rho \left( gh(z,t) + \frac{u^2(z,t)}{2} \right) \nonumber \\
		e_p  & = & \delta_p H = - \frac{(\star{\alpha^q}) 
(\star{\alpha^p})}{\rho} = -b(z)h(z,t)u(z,t)
\nonumber
\end{eqnarray}
which are the usual hydrodynamic pressure and the volume flow conjugated variables in the 1D hydraulic domain. Using these co-energy variables, the 1D SWE with non uniform width $b(z)$ reads:
\begin{equation}
	\label{eq:1DSWECovForm}
	\left[ \begin{array}{c} \dot{\alpha}^q \\  \dot{\alpha}^p \end{array} \right]  = \left[ \begin{array}{cc} 0 & \partial_z \\ - \partial_z & 0 \end{array} \right] \left[ \begin{array}{c}  e_q \\ e_p \end{array} \right] \,.
\end{equation}
All the geometry-dependent parameters have been embedded in the constitutive equations and  the structural interdomain coupling between potential and kinetic energies (i.e. between mass and momentum balance equations) in the SWE pops up again in the canonical 1D Stokes-Dirac structure in \eqref{eq:1DSWECovForm}.
\end{remark}

The Partitioned Finite Element Method directly applies to the  port-Hamiltonian system with variable coefficients. 
Here, the same procedure as in \S~\ref{subsec:wfSWE} is being followed. We begin with a weak formulation of (\ref{eq:1Dstrongformvariable}), then two complementary choices can be made.

If we choose to integrate by parts the {\em mass balance equation} only, i.e. the first line of the obtained weak form, we get  exactly the same finite-dimensional pHs as (\ref{eq:finitedim_waveeq1d}), but with $D$  a new $N_q \times N_p$ matrix, defined by 
 \begin{equation}
	\label{eq:DDefvariable}
		{D}  := \int_{[0,L]}{b(z)\, \frac{d \pmb{\phi}_{q}}{dz}(z) \, \pmb{\phi}^T_{p}(z) \dz} \,,
\end{equation}
and a new $N_q \times 2$ control matrix $B:=\left[\begin{matrix}b(0)\,\pmb{\phi}_q(0) & - b(L)\,\pmb{\phi}_q(L) \end{matrix} \right]$. The boundary control remains $\pmb{u}_{\partial}(t):= \left[\begin{matrix}e_p(0,t) & e_p(L,t) \end{matrix} \right]^T$, and the new collocated boundary observation reads  $\pmb{y}_{\partial}(t):= B^T \pmb{e}_q (t)= \left[\begin{matrix}b(0)\,e_q(0,t) & -b(L)\,e_q(L,t) \end{matrix} \right]^T$.

If instead, we choose to integrate by parts the {\em momentum balance equation} only, i.e. the second line of the obtained weak form, we get the following finite-dimensional pHs:
\begin{equation}
	\label{eq:finitedim_waveeq1dvariable2}
	\begin{split}
		{M_q} \dot{\pmb{q}}(t) & = {\widetilde{D}} \pmb{e}_p (t) \,, \\
		{M_p} \dot{\pmb{p}}(t)  & = {-\widetilde{D}^T} \pmb{e}_q (t) + \widetilde{B} \left[\begin{matrix}e_q(0,t) \\e_q(L,t) \end{matrix} \right] \,,
	\end{split}
\end{equation}
but with $\widetilde{D}$  another $N_q \times N_p$ matrix, defined by 
 \begin{equation}
	\label{eq:DDefvariable2}
	      \widetilde{D}  := - \int_{[0,L]}{ \pmb{\phi}_{q}(z)\, \frac{d}{dz}[b(z)\,\pmb{\phi}^T_{p}(z)] \dz} \,,
\end{equation}
and a new $N_p \times 2$ control matrix $\widetilde{B}:=\left[\begin{matrix}b(0)\,\pmb{\phi}_p(0) & - b(L)\,\pmb{\phi}_p(L) \end{matrix} \right]$.
The boundary control is now defined by $\pmb{u}_{\partial}(t):= \left[\begin{matrix}e_q(0,t) & e_q(L,t) \end{matrix} \right]^T$, and the new collocated boundary observation reads  $\pmb{y}_{\partial}(t):= \widetilde{B}^T \pmb{e}_p (t)= \left[\begin{matrix}b(0)\,e_p(0,t) & b(L)\,e_p(L,t) \end{matrix} \right]^T$.

Finally, note that in both the above cases, the following power balance is met:
$$
\dot{H_d}(t) :=\pmb{e}_q^T(t) {M_q} \dot{\pmb{q}}(t) + \pmb{e}_p^T(t) {M_p} \dot{\pmb{p}}(t)= \pmb{y}_{\partial}^T(t)\,\pmb{u}_{\partial}(t)\,.
$$

\color{black}
\subsection{Higher-order systems}
\label{ss-highorder}
In the previous sections, the PFEM was applied to first-order (1D and 2D) formally skew-symmetric differential operators. Indeed, the method seems to be much more general and can be applyied similarly to higher-order equations.

\paragraph{The Euler-Bernoulli beam equation can be written as a port-Hamiltonian system of second order}The equations are given by (see, e.g., \cite{CardosoRibeiro2016}):

\begin{equation}
		\label{eq:EulerBernoulliStrongform}
		\begin{split}
			\dot{x}_1(z,t) = - \frac{\partial^2}{\partial z^2} e_2(z,t) \,,\\
			\dot{x}_2(z,t) = \frac{\partial^2}{\partial z^2} e_1(z,t)\,,
		\end{split}
  \end{equation}
where $e_1$ and $e_2$ are obtained from the variational derivative of the Hamiltonian:
\begin{equation}
	H = \frac{1}{2}\int_0^L \left( x_1{}^2 + x_2{}^2  \right)\dz \,.
\end{equation}

From the definition of the variational derivatives, the time-derivative of the Hamiltonian is computed as:
\begin{align}
	\dot{H} =& \int_{z=0}^L \left(e_1 \dot{x_1} + e_2 \dot{x_2} \right) \dz\,, \nonumber\\
			=& \int_{z=0}^L\left(  - e_1 \partial_{z^2}^2e_2+ e_2 \, \partial_{z^2}^2 e_1 \right) \mathop{\mathrm{d}z}\,, \nonumber\\ \nonumber
			=& \int_{z=0}^L \left(\partial_z \left( - e_1 \, \partial_z (e_2) + \partial_z(e_1) \, e_2  \right)  \right)\dz \,, \nonumber\\
			= \nonumber & - e_1 (L,t) \, \partial_z (e_2)(L,t) + \partial_z(e_1) (L,t) e_2 (L,t) \\
			& + e_1 (0,t)\, \partial_z (e_2) (0,t) - \partial_z(e_1) (0,t) e_2 (0,t)   \,. \label{eq:Hdotinfdim}
\end{align}
Note that $\dot{H}$ depends only on the boundary values
of $e_1$ (vertical speed), $e_2$ (moment), $\partial_z e_1$ (rotation speed) and $\partial_z e_2$ (force). This motivates the definition of the boundary ports, which allows writing the infinite-dimensional equations as port-Hamiltonian systems.
From \eqref{eq:Hdotinfdim}, one possible definition for the boundary ports is as follows:
\begin{equation}
	\label{eq:boundary_ports_2ndorder}
		\pmb{y}_\partial:= \left[
		\begin{matrix} f_{1\partial}  \\ f_{2\partial}  \\ f_{3\partial} \\f_{4\partial} 
		\end{matrix} \right]
		:= \left[ \begin{matrix}  \partial_z \pmb{e}_1(L,t) \\  - \partial_z \pmb{e}_1(0,t)\\ -\pmb{e}_1(L,t)\\ \pmb{e}_1(0,t)
		\end{matrix} \right] \, , \,
		\pmb{u}_\partial = \left[
		\begin{matrix} e_{1\partial} \\  e_{2\partial} \\ e_{3\partial} \\ e_{4\partial}
		\end{matrix} \right]
		= \left[ \begin{matrix}  \pmb{e}_2(L,t) \\ \pmb{e}_2(0,t) \\ \partial_z \pmb{e}_2(L,t) \\  \partial_z \pmb{e}_2(0,t)
		\end{matrix} \right] \,.
\end{equation}

The final power balance ($\dot{H}$) can thus be written as:
\begin{align}
	\dot{H} = \pmb{y}_\partial^T \pmb{u}_\partial\,.
\end{align}

\paragraph{Weak-form representation of Euler-Bernoulli beam equation}
Let us use arbitrary test functions $v_1(z)$ and $v_2(z)$ and develop a weak form of \eqref{eq:EulerBernoulliStrongform}:
\begin{equation}
	\label{eq:EulerBernoulliWeakForm}
	\begin{split}
		\int_0^L v_1(z) \dot{x}_1(z,t) \dz = - \int_0^L v_1(z) \frac{\partial^2}{\partial z^2} e_2(z,t) \dz\,,\\
		\int_0^L v_2(z) \dot{x}_2(z,t) \dz = \int_0^L v_2(z) \frac{\partial^2}{\partial z^2} e_1(z,t) \dz\,,
	\end{split}
\end{equation}

Integrating the first equation by parts twice, we get the following partitioned weak form:
\begin{equation}
	\label{eq:EulerBernoulliWeakForm}
	\begin{split}
		\int_0^L v_1(z) \dot{x}_1(z,t) \dz   
		=&  - \int_0^L \frac{\partial^2}{\partial z^2}v_1(z) \,e_2(z,t) \dz 
		\\
		&		 + 
		 \left[
			 \begin{matrix}
					\partial_z v_1(L) & -\partial_z v_1(0) & -v_1(L) & v_1(0)
			 \end{matrix}
		 \right]
		 \left[
			\begin{matrix}
			   e_2(L,t) \\ e_2(0,t) \\ \partial_ze_2(L,t) \\ \partial_z e_2(0,t) 
			\end{matrix}
		\right]\,.\\
		\int_0^L v_2(z) \dot{x}_2(z,t) \dz =& \int_0^L v_2(z) \frac{\partial^2}{\partial z^2} e_1(z,t) \dz\,,
	\end{split}
\end{equation}

\paragraph{Finite-dimensional port-Hamiltonian system}
Similarly to the development in \S~\ref{subsec:wfSWE}, we chose finite-dimensional bases functions \eqref{eq:approx} as $\pmb{\phi_1}(z)$ and $\pmb{\phi_2}(z)$, for the variables with index 1 and 2, respectively.
From the substitution of the approximation functions in the weak form \eqref{eq:EulerBernoulliWeakForm}, we find:
\begin{equation}
	\label{eq:weakformFEM1DEB}
	\begin{split}
		{M_1} \dot{\pmb{x}}_1(t) = & \, - {D} \pmb{e}_2 (t) +  B 			
			\left[
			   \begin{matrix}
				  e_2(L,t) \\ e_2(0,t) \\ \partial_ze_2(L,t) \\ \partial_z e_2(0,t) 
			   \end{matrix}
		   \right] \,, \\
		   {M_2} \dot{\pmb{x}}_2(t) = & {D^T} \pmb{e}_1 (t)
	\end{split}
\end{equation}
 where $M_1$ and $M_2$ are square mass matrices (of size $N_1 \times N_1$ and $N_2 \times N_2$, respectively), equivalent to \eqref{eq:MqMpDef}.
The matrix $D$ is of size $N_1 \times N_2$:
\begin{equation}
	D := {\int_{z=0}^L{\left(\frac{\partial^2 \pmb{\phi}_{1}}{\partial z^2}(z) \right) \pmb{\phi}_{2}(z)^T \dz\,}} \,,
\end{equation}
and $B$ is an $N_1 \times 4$ matrix:
\begin{equation}
B := \left[\begin{matrix}  \frac{\partial \pmb{\phi}_1}{\partial z} (L) & - \frac{\partial \pmb{\phi}_1}{\partial z}(0) &-\pmb{\phi}_1(L) & \pmb{\phi}_1(0)  \end{matrix} \right]\,.
\end{equation}

Finally, the conjugated-output can also be written in terms of the previous $B$ matrix:

\begin{equation}
	\pmb{y}_\partial =
		\left[ \begin{matrix} \partial_z  \pmb{e}_1(L) \\  - \partial_z\pmb{e}_1(0)\\ -\pmb{e}_1(L)\\  \pmb{e}_1(0)
		\end{matrix} \right] = B^T \pmb{e}_1 \,.
\end{equation}

Defining the flow variables as $\pmb{f}_1(t) := -\dot{\pmb{x}}_1(t)$ and $\pmb{f}_1(t) := -\dot{\pmb{x}}_2(t)$, we find the following finite-dimensional Dirac structure representation:
\begin{equation}
	\label{eq:finitedimDirac_EB}
	\begin{split}
		\left[ \begin{matrix} M_1 & 0 \\ 0 & M_2 \end{matrix}\right] \left[ \begin{matrix} {\pmb{f}}_1(t) \\ {\pmb{f}}_2(t) \end{matrix}\right] & = \left[ \begin{matrix} 0 & -D \\ D^T & 0 \end{matrix}\right] \left[ \begin{matrix} {\pmb{e}}_q(t) \\ {\pmb{e}}_p(t) \end{matrix}\right] + \left[ \begin{matrix} -B \\ 0 \end{matrix}\right] \left[
			\begin{matrix}
			   e_2(L,t) \\ e_2(0,t) \\ \partial_ze_2(L,t) \\ \partial_z e_2(0,t) 
			\end{matrix}
		\right] \\
		\left[ \begin{matrix} \partial_z  \pmb{e}_1(L) \\  - \partial_z\pmb{e}_1(0)\\ -\pmb{e}_1(L)\\  \pmb{e}_1(0)
		\end{matrix} \right] & = \left[ \begin{matrix} B^T & 0 \end{matrix}\right] \left[ \begin{matrix} {\pmb{e}}_q(t) \\ {\pmb{e}}_p(t) \end{matrix}\right]
	\end{split}
\end{equation}

Following the same procedure presented in the previous sections for the 1D and 2D SWE, from the discretization of the Hamiltonian using the energy variables approximation spaces, one gets the underlying port-Hamiltonian dynamics for the approximated Euler-Bernoulli beam equations.

The analogue of the Euler-Bernoulli beam in 1D is the Kirchhoff plate in 2D, one can refer to \cite{BrugnoliAMM2019b} for the modelling as a port-Hamiltonian system using tensor calculus, and the application of PFEM to it, with various boundary controls ; note that the analogue of the Timoshenko beam in 1D is the Mindlin plate in 2D, and PFEM can also be applied to this model \cite{BrugnoliAMM2019a}.

\section{Conclusion and Open questions}
\label{sec:conclusions}

The Partitioned Finite Element Method provides a full-rank  structure-preserving representation of port-Hamiltonian systems in 2D and 3D: a general setting has been proposed here, written in the language of differential forms, and also translated into vector calculus for common PDE applications.
This method can be easily implemented thanks to ready to use FEM software to compute the matrices of the representation, which are all sparse.
It applies to  complex geometries, works in any  coordinate systems, and allows for space-varying coefficients; moreover higher order differential operators can also be tackled.
Although PFEM has already been succesfully applied to linear PDEs with quadratic Hamiltonian functionals, e.g. vibrating membranes and plates, here the methodology carries over to a non-linear PDE with non-quadratic and non-separable Hamiltonian functional, the irrotational Shallow Water Equation in 2D.

Future work will deal with mixed boundary control (possibly leading to differential algebraic problems as pHDAEs), and mathematical convergence analysis (choice of the finite element bases and theoretical rate of convergence). Some worked-out 2D test cases will be studied on coupled systems, e.g. fluid-structure interaction (FSI), or thermal-structure coupling. Lastly, structure-preserving model reduction techniques will be tested on the high-fidelity finite-dimensional systems obtained by PFEM, see e.g. \cite{EggKugLilMarMeh18} for pHs or \cite{Mehrmann2019} for pHDAEs. The reduced order system will then be most useful to apply dedicated control laws for pHs, like IDA-PBC, which do take advantage of the specific structure of these  dynamical systems with collocated inputs and outputs, see e.g. \cite{Ortega2008}.

\vspace{1cm}
\noindent
    {\bf THANKS} 
    Part of this work has been performed in the framework of the Collaborative Research DFG and ANR project INFIDHEM n$^\circ$ ANR-16-CE92-0028 ({\tt http://websites.isae.fr/infidhem}).

\vspace{1cm}
\begin{center}
  \rule{4cm}{0.5pt}
\end{center}



\end{document}